\documentclass[12pt,a4paper]{article}

\usepackage[T1]{fontenc}
\usepackage[latin1]{inputenc}
\usepackage{graphicx}
\usepackage{epsfig}
\usepackage{fancyhdr,fancybox}
\usepackage{amsfonts}
\usepackage{amsmath}
\usepackage{amssymb}
\usepackage{verbatim}
\usepackage{subfigure}
\usepackage[francais,english]{babel}
\usepackage[latin1]{inputenc}
\usepackage{latexsym}
\usepackage{multicol}

\usepackage{graphicx}
\usepackage{subfigure}

\newcommand{\R}{\mathbb{R}}

\newcommand{\inter}{{\rm int}\;}
\newcommand{\N}{\mathbb{N}}
\newcommand{\E}{\mathbb{E}}

\def\acknowledgement{\par\addvspace{17pt}\small\rmfamily
\trivlist\if!\ackname!\item[]\else
\item[\hskip\labelsep
{\it\ackname}]\fi}
\def\endacknowledgement{\endtrivlist\addvspace{6pt}}

\newcommand\ackname{Acknowledgements\runinend}
\def\runinend{.}

\newcommand{\kommentar}[1]{}
\newcommand{\ind}{{1\hspace{-1mm}{\rm I}}}

\newtheorem{thm}{Theorem}[section]

\newtheorem{lem}[thm]{Lemma}
\newtheorem{cor}[thm]{Corollary}

\begin{document}
\pagestyle{fancy}
\title{Scaling limits for shortest path lengths along the edges of stationary tessellations\\[0.2cm] Supplementary material}
\author{Florian Voss$^{1}$ \and Catherine Gloaguen $^2$ \and Volker Schmidt$^1$}
\date{\today}
\maketitle


\bigskip

\begin{abstract}
We consider spatial stochastic models, which can be applied e.g.
to telecommunication networks with two hierarchy levels. In
particular, we consider Cox processes $X_L$ and $X_H$ concentrated
on the edge set $T^{(1)}$ of a random tessellation $T$, where the
points $X_{L,n}$ and $X_{H,n}$ of $X_L$  and $X_H$ can describe
the locations of low--level and high--level network components,
respectively, and $T^{(1)}$ the underlying
infrastructure of the network, like road systems,
railways, etc. Furthermore, each point $X_{L,n}$ of $X_L$ is
marked with the shortest path along the edges of $T$ to the
nearest (in the Euclidean sense) point of $X_H$. We investigate
the typical shortest path length $C^*$ of the resulting marked
point process, which is an important characteristic e.g. in
performance analysis and planning of telecommunication networks. In particular, we
 show that the distribution of $C^*$ converges to simple parametric
 limit distributions if a scaling factor $\kappa$ converges to
zero and infinity, respectively.  This can be used to approximate the density of $C^*$
by analytical formulae for a wide range of $\kappa$.\\[0.5cm]
\vspace{0.05cm}
\begin{tabbing}
\noindent {\it AMS {\rm 2000} subject classification}$\,$: {\sc
60D05, 60G55, 60F99, 90B1}\\
\end{tabbing}
\end{abstract}

\vspace*{-1cm}
\stepcounter{footnote} \footnotetext{Institute of Stochastics, Ulm University, 89069 Ulm, Germany\\
\hspace*{0.9cm}e-mail: florian.voss@uni-ulm.de, volker.schmidt@uni-ulm.de}
\stepcounter{footnote} \footnotetext{Orange Labs, 92131 Issy les Moulineaux Cedex 9, France\\
\hspace*{0.9cm}e-mail: catherine.gloaguen@orange-ftgroup.com}


\numberwithin{equation}{section}
\newpage

\section{Introduction}\label{sec.two.one}

Asymptotic properties of spatial stochastic models are considered, which can be applied e.g.
in the analysis and planning of telecommunication networks.
More precisely, we consider stochastic models for networks with two
hierarchy levels, i.e., there are network components of two different kinds:
low-level components (LLC) and
high-level components (HLC). The locations of both HLC and LLC are represented
by points in the Euclidean plane $\R^2$. We then associate with each HLC a certain subset of $\R^2$ which is called its serving zone.
This is done in such a way that the serving zones of the HLC are
disjoint convex polygons which cover the whole $\R^2$. Each LLC is
linked to the HLC in whose serving zone the LLC is
located. In particular, we assume that the serving
zones are constructed as the cells of the Voronoi tessellation
with respect to the locations of HLC. This is equivalent to link each LLC to its
nearest HLC, where "nearest" means with respect to the Euclidean distance.
Furthermore, we  assume that the
HLC and LLC are located on the edges of a random geometric graph,
where the link from a LLC to its nearest HLC is
assumed to be the shortest path along the edges of that graph.
In the case of telecommunication networks the edges of the random geometric
graph represent the underlying infrastructure, e.g.
an inner-city street system.

Thus, we study a class of stochastic network models which has been
introduced in~\cite{Glo02} as
the Stochastic Subscriber Line Model (SSLM) for urban access networks.
 Note that the SSLM is a model from
stochastic geometry which provides tools for the description of
geometric features of the network. Based on this model, stochastic
econometrical analysis can be done for real telecommunication networks, e.g.
connection costs for access networks can be determined,
see~\cite{Glo07,Glo09,Voss08b,Voss09a}, where we focus on the case
that the infrastructure of the network is modeled by the edge set of a stationary random
tessellation and both the HLC and LLC are modeled by Cox
processes concentrated on this edge set.
Then we are especially interested in the shortest path length
along the edge set between LLC and HLC, which is an important
performance characteristic in
cost and risk analysis as well as in strategic planning  of wired telecommunication. In order to
define an appropriately chosen (global) distribution of the
shortest path length we regard the so-called typical shortest path
length $C^*$. It can be interpreted as the length of
the shortest path from a location of LLC, which is chosen at
random among all locations of LLC, and its nearest HLC. 
We are
then interested in the asymptotic behaviour of the distribution of $C^*$ for two extreme
cases of model parameters. In particular, we  
 show that the distribution of $C^*$ converges to simple parametric
 limit distributions if a scaling factor $\kappa$ converges to
zero and infinity, respectively.  This can be used to approximate the density of $C^*$
by analytical formulae for a wide range of $\kappa$
which is a great advantage e.g.
for the econometrical analysis of real telecommunication networks, see \cite{Glo09}.
The mathematical techniques, which we
exploit in order to derive our main results presented in Theorems~\ref{the.gam.zer} and \ref{the.gam.inf},
include Palm calculus and Poisson approximation for stationary point processes, Kingman's subadditive ergodic theorem,
and the  generalized Blaschke-Petkantschin formula from geometric measure theory.

The paper is organized as follows. In Section \ref{sec.modelling} we give a short description of the particular
stochastic network model considered in the present paper.
 Then, in Section \ref{sec.asy.sho}, we present the main results stated in Theorems~\ref{the.gam.zer} and \ref{the.gam.inf}.
The proof of Theorem~\ref{the.gam.inf} is given in Section~\ref{sec.fou}, where some details are postponed to the Appendix.
In Section~\ref{sec.sec.fiv}, it is shown that the mixing and integrability
conditions of  Theorems~\ref{the.gam.zer} and \ref{the.gam.inf} are fulfilled for various examples of random tessellations.
  Some extensions of our results to other performance characteristics, more general classes of random geometric graphs,
  and more general connection rules are discussed in Section~\ref{sub.fur.tes}.
  Finally, Section~\ref{sec.con.out} concludes the paper and gives an outlook to possible future research.


\section{Stochastic modelling of hierarchical networks}
\label{sec.modelling}\setcounter{equation}{0}

To begin with we give a short description of the particular
stochastic network model considered in the present paper.
For more details on this model see also \cite{Glo07}.
Moreover, we briefly explain the mathematical background and
introduce the notation we are using. For further details on spatial point
processes and random tessellations, see e.g.
\cite{Daley0307,Okabe00,ScW08,Stoyan95}. Surveys on applications of tools from stochastic geometry
to spatial stochastic modelling of telecommunication networks can be found e.g. in \cite{Hae09,Zuy09}.

\subsection{Marked point processes}

First we recall some basic notions and results regarding marked
point processes in $\R^2$. They can be used to model locations of customers
or equipments in telecommunication networks. Let $\mathcal{B}^2$
denote the family of Borel sets of $\R^2$ and $N$ the family of
all simple and locally finite counting measures on
$\mathcal{B}^2$. Note that each $\varphi\in N$ can be represented
by the sequence $\{x_n\}$ of its atoms, i.e. $\varphi=\sum_n
\delta_{x_n}$, where $\delta_x$ is the Dirac measure with
$\delta_x(B)=1$ if $x\in B$ and $\delta_x(B)=0$ if $x\not\in B$.
Let $\mathcal{N}$ denote the $\sigma$-algebra of subsets of $N$ generated by the
sets $\{\varphi\in N: \varphi(B)=j\}$ for $j\in\N$ and
$B\in\mathcal{B}^2$. The shift operator $t_x:N\mapsto N$ is
defined by $t_x\varphi(B)=\varphi(B+x)$ for $x\in\R^2$ and
$B\in\mathcal{B}^2$, where $B+x=\{x+y: y\in B\}$. Then a point
process $X$ is a random element of the measurable space
$(N,\mathcal{N})$, where we identify $X$ with the sequence
$\{X_n\}$ of its (random) atoms, writing $X=\{X_n\}$ for brevity.

Let $\mathbb{M}$ be a Polish space with its Borel $\sigma$-algebra
$\mathcal{B}_\mathbb{M}$. Then we use the notation $N_\mathbb{M}$
for the family of all counting measures on $\mathcal{B}^2\otimes
\mathcal{B}_\mathbb{M}$ which are simple and locally finite in the
first component. Note that the atoms $(x_n,m_n)$ of the counting
measure $\psi=\sum_n \delta_{(x_n,m_n)}\in N_\mathbb{M}$ have two
components: the location $x_n\in\R^2$ and the mark
$m_n\in\mathbb{M}$. The $\sigma$-algebra $\mathcal{N}_\mathbb{M}$
is defined in the same way as above and the shift operator
$t_x:N_\mathbb{M}\mapsto N_\mathbb{M}$ translates the first
component of the atoms of $\psi\in N_\mathbb{M}$ by $-x$, i.e.
$t_x(\psi)=\sum_n \delta_{(x_n-x,m_n)}$. A random element
$X=\{(X_n,M_n)\}$ of $(N_\mathbb{M},\mathcal{N}_\mathbb{M})$ is
then called a marked point process.

\subsection{Palm distributions}
\label{subsec.pal.dis}

Stationarity and isotropy of (marked) point processes
are defined in the usual way, i.e., assuming the invariance of their distributions with respect
to arbitrary translations and rotations around the origin, respectively. By $\lambda >0$ we denote the
intensity of a stationary marked point process $X=\{(X_n,M_n)\}$, i.e. $\lambda=\E\#\{n:\,X_n\in[0,1]^2\}$,
and the Palm mark distribution
$\mathbb{P}^{o}_X:\mathcal{B}_\mathbb{M}\to[0,1]$ of $X$ is given
by
\begin{equation}\label{pal.mar.dis}
    \mathbb{P}^{o}_X(G)=\frac{\E\#\{n:X_n\in[0,1)^2,M_n\in G\}}{\lambda}\;,\qquad G\in\mathcal{B}_\mathbb{M}\,.
\end{equation}
A random variable $M^*$ distributed according to
$\mathbb{P}^{o}_X$ is called the typical mark of~$X$.

Furthermore, two jointly stationary marked point processes
$X^{(1)}=\{(X_n^{(1)},M_n^{(1)})\}$ and
$X^{(2)}=\{(X_n^{(2)},M_n^{(2)})\}$ with intensities $\lambda_{1}$
and $\lambda_{2}$ and mark spaces $\mathbb{M}_1$ and
$\mathbb{M}_2$, respectively, will be considered as random element
$Y=(X^{(1)},X^{(2)})$ of the product space
$N_{\mathbb{M}_1,\mathbb{M}_2}=N_{\mathbb{M}_1}\times
N_{\mathbb{M}_2}$. The Palm distribution
$\mathbb{P}^{*}_{X^{(i)}}$ of $Y$ with respect to the $i$-th
component, $i=1,2$, is then defined on
$\mathcal{N}_{\mathbb{M}_1}\otimes
\mathcal{N}_{\mathbb{M}_2}\otimes \mathcal{B}_{\mathbb{M}_i}$ by
\begin{equation}\label{def.pal.dis}
    \mathbb{P}^{*}_{X^{(i)}}(A\times G)=\frac{\E\#\{n:X_n^{(i)}\in[0,1)^2,M_n^{(i)}\in G, t_{X_n^{(i)}}Y\in A\}}{\lambda_i}\,,
\end{equation}
where $A\in \mathcal{N}_{\mathbb{M}_1}\otimes
\mathcal{N}_{\mathbb{M}_2}$ and $G\in \mathcal{B}_{\mathbb{M}_i}$.
Note that the Palm mark distribution $\mathbb{P}^{o}_{X^{(i)}}$ of
$X^{(i)}$ can be obtained from $\mathbb{P}^{*}_{X^{(i)}}$ as a
marginal distribution.

\kommentar{%
Finally, we need Neveu's exchange formula
for jointly stationary marked point processes. Using the notation
introduced above, and $\psi=(\psi^{(1)},\psi^{(2)})$ for the
elements of $N_{\mathbb{M}_1,\mathbb{M}_2}$, this formula takes
the following form (see e.g. \cite{Mai04,NE76}).
\begin{lem}
\label{th.nev.exc}  For any measurable $f:\R^2\times
\mathbb{M}_1\times \mathbb{M}_2\times
N_{\mathbb{M}_1,\mathbb{M}_2}\rightarrow [0,\infty)$, it holds
that
\begin{equation}
\begin{split}
\label{eq.nev.exc} \lambda_{1} &
\int_{N_{\mathbb{M}_1,\mathbb{M}_2}\times
\mathbb{M}_1}\,\int_{\R^2\times \mathbb{M}_2}
f(x,m_1,m_2,t_x\psi)\,\,\psi^{(2)}(d(x,m_2))\,\mathbb{P}^{*}_{X^{(1)}}(d(\psi, m_1))\\
& = \lambda_{2} \int_{ N_{\mathbb{M}_1,\mathbb{M}_2}\times
\mathbb{M}_2}\,\int_{\R^2\times \mathbb{M}_1}
f(-x,m_1,m_2,\psi)\,\,
\psi^{(1)}(d(x,m_1))\,\mathbb{P}^{*}_{X^{(2)}}(d(\psi,m_2))\,.
\end{split}
\end{equation}\label{nev.exc.for}
\end{lem}
}%

\subsection{Random tessellations}
\label{subsec.models}

As a model for the underlying random geometric graph we
consider the edge set of random tessellations of $\R^2$. Note that
a random tessellation $T$ is a locally finite partition
$\{\Xi_n\}$ of $\R^2$ into random (compact and convex) polygons
$\Xi_n$, which are called the cells of $T$. We can also regard $T$
as a marked point process $\{(\alpha(\Xi_n),\Xi_n^o)\}$, where the
shifted cells $\Xi^o_n=\Xi_n-\alpha(\Xi_n)$ contain the origin.
The points $\alpha(\Xi_n)\in\Xi_n\subset\R^2$ are then called the
nuclei of the cells $\Xi_n$ of $T$.  Furthermore, we can identify
$T$ with the edge set $T^{(1)}=\bigcup_n\partial\Xi_n$ of $T$.
Note that $T^{(1)}$ is a random closed set in $\R^2$, i.e.,
$T^{(1)}$ is a random element of
$(\mathcal{F},\mathcal{B}(\mathcal{F}))$, where $\mathcal{F}$
denotes the family of all closed subsets of $\R^2$ and
$\mathcal{B}(\mathcal{F})$  is the smallest $\sigma$-algebra of
subsets of $\mathcal{F}$ which contains the ,,hitting sets''
$\mathcal{F}_C=\{B\in \mathcal{F}:\,B\cap C\not=\emptyset\}$ for
all compact $C\in\mathcal{B}^2$.

If $T$ is stationary, i.e., $T^{(1)}\stackrel{\rm d}{=}T^{(1)}+x$ for each $x\in \R^2$, then the intensity $\gamma$ of $T$ is
defined as $\gamma=\E \nu_1(T^{(1)}\cap [0,1]^2)$, i.e. the
mean length of $T^{(1)}$ per unit area, where
 $\nu_1$ denotes the 1-dimensional Hausdorff measure.
In the following we always assume that $T$ is a (normalized)
stationary tessellation with $\E\nu_1(T^{(1)}\cap [0,1]^2)=1$.
Furthermore, for each $\gamma>0$ we consider the scaled
tessellation $T_\gamma$ with intensity $\gamma$ which is defined
by $T_\gamma= T/\gamma$, i.e., we scale the edge set $T^{(1)}$
with $1/\gamma$ getting $T_\gamma^{(1)}$ such that
$\E\nu_1(T^{(1)}_\gamma\cap [0,1]^2)=\gamma$.

A random tessellation  $T$ is called isotropic if the distribution of $T^{(1)}$ is invariant
with respect to rotations around the origin. Furthermore, a stationary tessellation $T$ is called mixing if
\begin{displaymath}
    \lim_{|x|\longrightarrow \infty}\mathbb{P}(T^{(1)}\in A, T^{(1)}+x\in A^\prime)=\mathbb{P}(T^{(1)}\in A)\;\mathbb{P}(T^{(1)}\in A^\prime)
\end{displaymath}
for any     $A,A^\prime\in\mathcal{B}(\mathcal{F})$.
Note that for any $T$ which is mixing it holds that
\begin{equation}\label{def.erg.tes}
\mathbb{P}(T^{(1)}\in A)=1\quad\mbox{or}\quad \mathbb{P}(T^{(1)}\in A)=0\qquad\mbox{for each $A\in\mathcal{I}(\mathcal{F})$,}
\end{equation}
 where $\mathcal{I}(\mathcal{F})$ denotes the sub-$\sigma$-algebra of invariant sets
of $\mathcal{B}(\mathcal{F})$, i.e. $A+x=A$ for all $A\in\mathcal{I}(\mathcal{F})$ and $x\in\R^2$.
A stationary tessellation $T$ which satisfies condition (\ref{def.erg.tes}) is said to be ergodic.

\subsection{Cox processes on edge sets}
\label{subsec.low.upp}

For any $\gamma>0$, we consider Cox point processes
$X_H=\{X_{H,n}\}$ and $X_L=\{X_{L,n}\}$ concentrated on
$T^{(1)}_\gamma$, in order to model the locations of
 HLC and LLC, respectively. In particular, we assume that
$X_H$ is a Cox process on $T^{(1)}_\gamma$ with linear intensity
$\lambda_\ell$ which is constructed by placing homogeneous Poisson processes
on the edges of $T_\gamma$ with linear intensity $\lambda_\ell$.
The random driving measure
$\Lambda_{X_H}:\mathcal{B}^2\longrightarrow [0,\infty]$ of $X_H$
is then given by
\begin{equation}\label{ran.dri.mea}
    \Lambda_{X_H}(B)=\lambda_\ell\nu_1(B\cap T_\gamma^{(1)}),\qquad B\in \mathcal{B}^2\,.
\end{equation}
Analogously, $X_L$ is a Cox process on $T^{(1)}_\gamma$ with
linear intensity $\lambda^\prime_\ell$ which is constructed in the
same way, i.e., by placing Poisson processes on the edges of
$T_\gamma$ with linear intensity $\lambda^\prime_\ell$. Thus,
$X_H$ and $X_L$ are Cox processes concentrated on the same edge
set $T^{(1)}_\gamma$, where we assume that $X_H$ and $X_L$ are
conditionally independent given $T_\gamma$. Furthermore, note that
$X_H$ and $X_L$ are stationary, isotropic, and ergodic if $T$ is
stationary, isotropic, and ergodic, respectively. The planar
intensities $\lambda$ and $\lambda^\prime$ of $X_H$ and $X_L$ are
given by $\lambda=\lambda_\ell\gamma$ and
$\lambda^\prime=\lambda_\ell^\prime\gamma$.

\subsection{Serving zones and shortest paths}
\label{subsec.ser.zon}

Let $T_H=\{\Xi_{H,n}\}$ denote the Voronoi tessellation induced by
the points $X_{H,n}$ of the Cox process $X_H=\{X_{H,n}\}$, i.e.
\[
\Xi_{H,n}=\{x\in\R^2:\,|x-X_{H,n}|\le |x-X_{H,m}| \,\mbox{for all $m\not= n$}\}\,,
\]
where $|\,\cdot\,|$ denotes the Euclidean norm.
The
cells $\Xi_{H,n}$ of $T_H$  are considered to be the serving zones of
HLC. By means of the four modelling components $T_\gamma$, $X_H$,
$X_L$ and $T_H$ we can construct the marked point process
$X_{L,C}=\{(X_{L,n},C_n)\}$, where the mark $C_n$ is the length of the
shortest path from $X_{L,n}$ to $X_{H,j}$ along the edge set
$T^{(1)}_\gamma$ of $T_\gamma$ provided that $X_{L,n}\in
\Xi_{H,j}$.

Thus, each LLC is connected to its nearest HLC in the Euclidean sense and not in the shortest path sense. However, 
for applications this is a reasonable assumption since the planning of telecommunication networks is complicated and existing 
networks have evolved for long periods. 
Therefore, it is unrealistic to assume that serving zones are defined with respect to the shortest path distance and
it is appropriate to use a simpler rule. Furthermore,
analysis of real data has shown that the approach considered in the present paper is realistic (\cite{Glo09}).

\begin{figure}
\begin{center}
\subfigure[PVT as infrastructure model]{\includegraphics[width = 7cm, viewport = 0 0 116 126, clip] {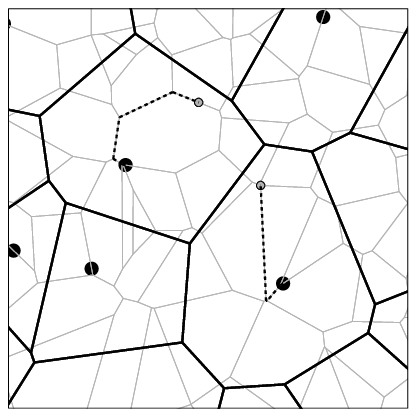}}
\subfigure[PLT as infrastructure model]{\includegraphics[width = 7cm, viewport = 42 75 158 201, clip] {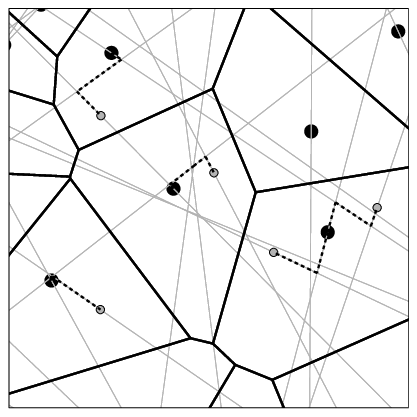}}
{\caption {Higher-level components with their serving zones (black)
and lower-level components (grey with black boundary) with shortest paths (dashed)
along the edge set (grey).} \label{fig.sp}}
\end{center}
\end{figure}

It is not difficult to show that $X_{L,C}$ is a
stationary and isotropic marked point process if $T_\gamma$ is
stationary and isotropic, respectively. Realizations of
service zones and shortest paths are displayed in
Figure~\ref{fig.sp}(a) and (b) for $T_\gamma$ being a Poisson-Voronoi tessellation (PVT) and a Poisson line tessellation (PLT),
respectively.

The model characteristic we are mainly interested in is the
distribution of the typical mark $C^*$ of $X_{L,C}$. Thus, we are
interested in the Palm mark distribution $\mathbb{P}^{o}_{X_{L,C}}$ of
$X_{L,C}$, i.e., the distribution of the typical shortest path length.

Note that the realizations of $X_{L,C}$ can be constructed from the
corresponding realiza\-tions of $X_L$ and $X_{H,S}$, where
$X_{H,S}=\{(X_{H,n},S^o_{H,n})\}$ is a stationary marked point process
with marks $S_{H,n}^o=(T_\gamma^{(1)}\cap\Xi_{H,n})-X_{H,n}$.
Thus, instead of $X_{L,C}$, we can consider the vector $Y=(X_L, X_{H,S})$
and the Palm distribution $\mathbb{P}_{X_L}^*$ of $Y$ with respect
to $X_L$, which has been introduced in (\ref{def.pal.dis}).
  Let $(X_L^*,\widetilde{X}_{H,S})$
be distributed according to $\mathbb{P}_{X_L}^*$, where we use the
notation
$\widetilde{X}_{H,S}=\{(\widetilde{X}_{H,n},\widetilde{S}_{H,n}^o)\}$
and
\begin{equation}\label{def.teh.gam}
\widetilde{T}_\gamma^{(1)}=\bigcup_{n\ge
1}\bigl(\widetilde{S}_{H,n}^o+\widetilde{X}_{H,n}\bigr)\,.
\end{equation}
Note that $\widetilde{X}_H=\{\widetilde{X}_{H,n}\}$ is a
(non-stationary) Cox process on $\widetilde{T}^{(1)}_\gamma$ with
linear intensity $\lambda_\ell$. Moreover, by
$\widetilde{X}_{H,0}$ we denote the closest point (in the
Euclidean sense) of $\{\widetilde{X}_{H,n}\}$ to the origin. Then,
the typical shortest path length $C^*$ can be given by
$C^*=c(\widetilde{X}_{H,0})$, where $c(\widetilde{X}_{H,0})$
denotes the length of shortest path from the origin to
$\widetilde{X}_{H,0}$, along the edges of
$\widetilde{T}_\gamma^{(1)}$. In the following we always assume that
the joint distribution of $C^*, \widetilde{X}_H$ and $\widetilde{T}_\gamma$ is
given by $\mathbb{P}^*_{X_L}$.



\section{Limit theorems for the typical shortest path length}
\label{sec.asy.sho}

We investigate the asymptotic behavior of the distribution of $C^*$
for two different cases: $\gamma\to 0$ with $\lambda_\ell$  fixed and $\gamma/\lambda_\ell\to\infty$
with $\gamma\lambda_\ell$ fixed,
i.e., unboundedly sparse edge sets and unboundedly dense edge
sets, respectively. For $\gamma\to 0$, we show in
Theorem~\ref{the.gam.zer}
that the distribution of $C^*$ converges weakly to an exponential distribution, where no
specific assumption on the underlying stationary tessellation $T$
is needed. Furthermore, for $\gamma\to\infty$ and $T$ being a
stationary and isotropic random tessellation which is mixing, we
get in Theorem~\ref{the.gam.inf} that the distribution of $C^*$
converges weakly to a Weibull distribution.

\subsection{Scaling invariance property}\label{sebsec.aux.con}

Recall that the stochastic network model introduced in
Section~\ref{sec.modelling} and, in particular, the distribution
of $C^*$ is fully specified by $T$, $\gamma$,
$\lambda_\ell$ and $\lambda_\ell^\prime$. Moreover,
it can be shown (see e.g. \cite{Glo07,Voss08b}) that the
distribution of $C^*$ does not depend on $\lambda_\ell^\prime$.
Therefore,  we only regard the parameters $\gamma$ and
$\lambda_\ell$ in the following.  Sometimes we use the notation
$C^*=C^*(\gamma,\lambda_\ell)$ to emphasize that the distribution
of $C^*$ depends on $\gamma$ and $\lambda_\ell$.

\begin{figure}
\centering
\subfigure[$\kappa$ = 0.5]{
\includegraphics[width=7.0cm]{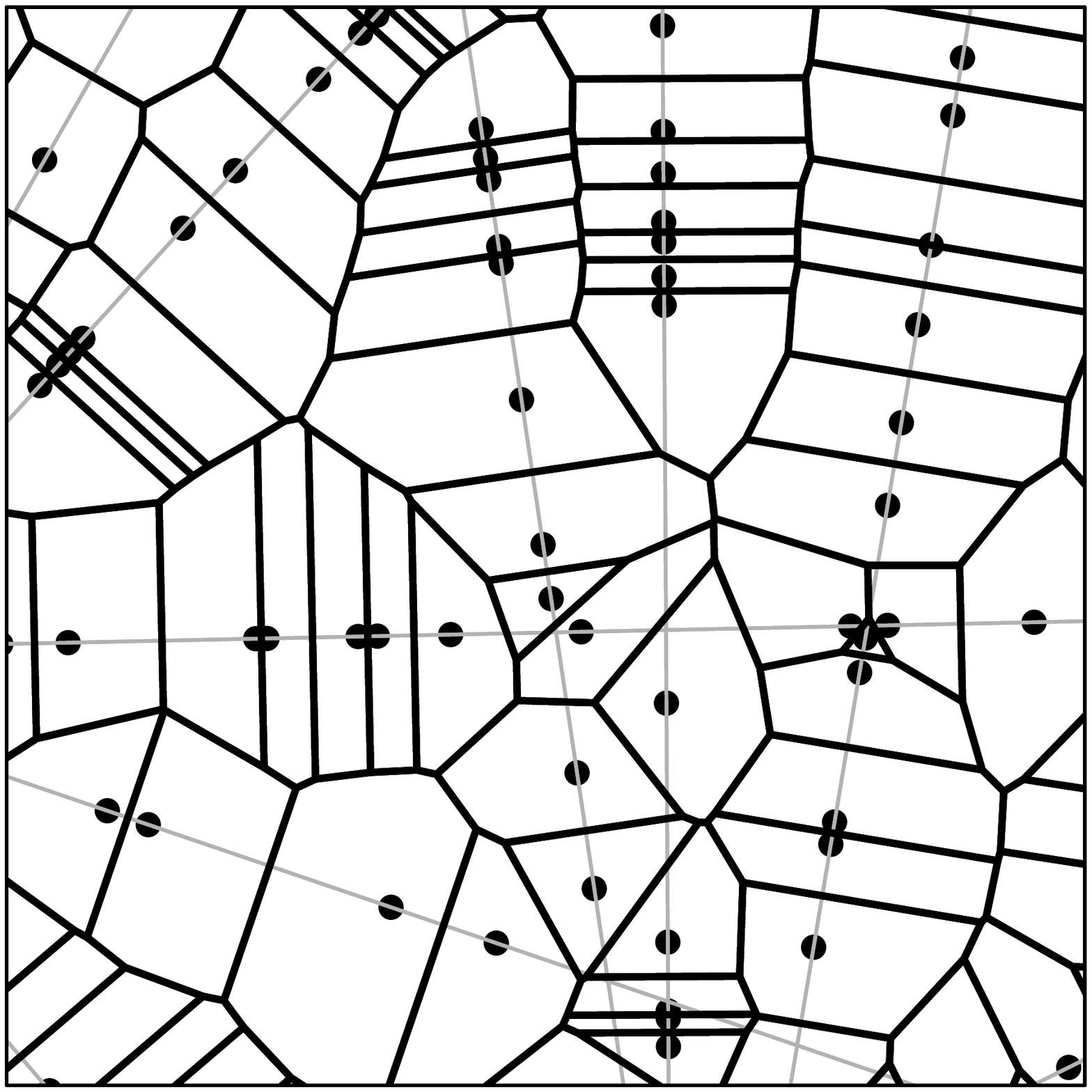}}
\subfigure[$\kappa$ = 1000]{
\includegraphics[width=7.0cm]{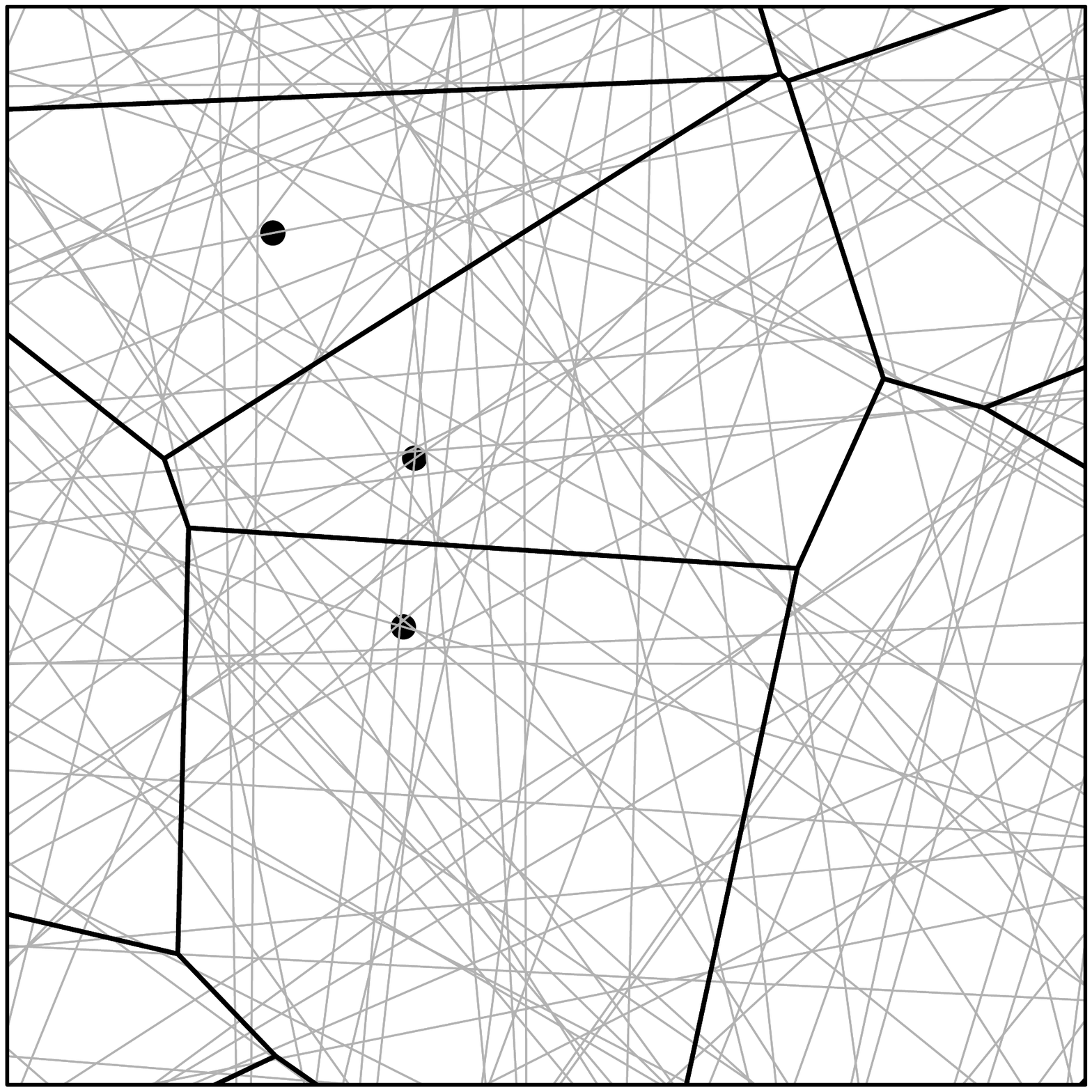}
} \caption{Realizations of $X_{H,S}=\{(X_{H,n},S^o_{H,n})\}$ for
extremal values of $\kappa$} \label{fig:conv}
\end{figure}

Furthermore, a scaling invariance property holds for this model.
If the value of the quotient $\kappa=\gamma/\lambda_\ell$ is
constant, then the structure of $X_{H,S}$ is fixed, but on
different scales for different parameter vectors
$(\gamma,\lambda_\ell)=(\kappa\lambda_\ell,\lambda_\ell)$.
We are interested in the limiting behavior of the distribution of
$C^*$ for $\kappa\to 0$ with $\lambda_\ell$ fixed and
for $\kappa\to \infty$ with
$\lambda=\lambda_\ell\gamma$ fixed.
 In Figure \ref{fig:conv}  realizations of $X_{H,S}$ are shown for two
(extremely small and large) values of $\kappa$, where the
realization of $T$ is sampled from a PLT. One can see that for
small values of $\kappa$ the segment systems within the serving
zones mainly consist of one single segment only, whereas for large
values of $\kappa$ the networks inside the serving zones become
rather dense.

\subsection{Asymptotic exponential distribution for $\kappa\to 0$}

First we regard the case that $\kappa=\gamma/\lambda_\ell\to 0$
with $\lambda_\ell$ fixed, i.e., $\gamma\to 0$.
\begin{thm}\label{the.gam.zer}
    Let $T$ be an arbitrary stationary tessellation. Then, for any
    fixed $\lambda_\ell>~0$, it holds that
    \begin{equation}
        C^*(\gamma,\lambda_\ell)\stackrel{\rm d}{\to } Z\qquad\mbox{as}\quad \gamma\to 0\,,
    \end{equation}
    where $\stackrel{\rm d}{\to}$ denotes convergence in distribution and $Z\sim {\rm Exp}(2\lambda_\ell)$, i.e.,
    the random variable $Z$ is exponentially distributed
    with expectation $(2\lambda_\ell)^{-1}$.
\end{thm}
\begin{proof} Let
$R_\gamma= \max\{r>0:\, B(o,r)\cap
\widetilde{L}_\gamma^o=B(o,r)\cap \widetilde{T}_\gamma^{(1)}\}$,
where $B(o,r)$ denotes the ball centered at the origin with radius
$r$ and $\widetilde{L}_\gamma^o$ is the segment containing the
origin of the random edge set $\widetilde{T}^{(1)}_\gamma$
introduced in (\ref{def.teh.gam}). It is not difficult to see that
\begin{equation}\label{con.err.inf}
\lim\limits_{\gamma\to
0}R_\gamma=\infty\qquad\mbox{ a.s.}
\end{equation}
Recall that $C^*=c(\widetilde{X}_{H,0})$, where
$\widetilde{X}_{H,0}$ is the closest point  to the origin of the
point process $\widetilde{X}_H=\{\widetilde{X}_{H,n}\}$ of HLC
under $\mathbb{P}^*_{X_L}$, and note that the values of the
distribution function $F_{C^*}:(0,\infty)\to(0,1)$ of $C^*$ can be
written as
\begin{eqnarray*}
    F_{C^*}(x)&=&\mathbb{P}(\widetilde{X}_{H,0}\in B(o,R_\gamma))\;
    \mathbb{P}(C^*\le x \mid \widetilde{X}_{H,0}\in B(o,R_\gamma))\\
            &&+\;\mathbb{P}(\widetilde{X}_{H,0}\not\in B(o,R_\gamma))\;
            \mathbb{P}(C^*\le x \mid \widetilde{X}_{H,0}\not\in B(o,R_\gamma))
\end{eqnarray*}
for each $x\ge 0$. It can be shown (see e.g. \cite{FL08}) that
$\widetilde{X}_H$ is a Cox process which is a homogeneous
Poisson process with linear intensity $\lambda_\ell$ on the
edges of $\widetilde{T}^{(1)}_\gamma$. This implies that
\[
 \mathbb{P}(C^*\le x \mid \widetilde{X}_{H,0}\in B(o,R_\gamma))=
    \frac{\mathbb{P}(\min\{Z_1,Z_2\}\le x,\min\{Z_1,Z_2\}\le R_\gamma )}
    {\mathbb{P}(\widetilde{X}_{H,0}\in B(o,R_\gamma))}
\]
for each $x>0$, where the random variables $Z_1$ and $Z_2$ are independent,
exponentially distributed with parameter $\lambda_\ell$ and independent of $R_\gamma$.
Furthermore, we get that
\[
\mathbb{P}(\widetilde{X}_{H,0}\not\in
B(o,R_\gamma))=\mathbb{P}(\min\{Z_1,Z_2\}> R_\gamma)
=\mathbb{E}\exp(-2\lambda_\ell R_\gamma)\,,
\]
since $\min\{Z_1,Z_2\}$ is exponentially distributed with
parameter $2\lambda_\ell$ and independent of $R_\gamma$. Thus,
using (\ref{con.err.inf}), it follows that
\[
\lim\limits_{\gamma\to
0}\mathbb{P}(\widetilde{X}_{H,0}\not\in B(o,R_\gamma))=
0\qquad\mbox{and}\qquad \lim\limits_{\gamma\to 0}
\mathbb{P}(\widetilde{X}_{H,0}\in B(o,R_\gamma))= 1
\]
 and, consequently, $\lim_{\gamma\to 0}
F_{C^*}(x)=\mathbb{P}(\min\{Z_1,Z_2\}\le
x)\;=\;1-\exp(-2\lambda_\ell x)$  for each $x\ge 0$.\end{proof}

Note that the case $\kappa=\gamma/\lambda_\ell\to 0$ with $\gamma$ fixed and $\lambda_\ell\to \infty$
can be treated in the following way. Due to the scaling invariance property mentioned
in Section~\ref{sebsec.aux.con} we have
\[
\lambda_\ell\;\, C^*(\gamma,\lambda_\ell)\stackrel{\rm d}{=} C^*(\gamma/\lambda_\ell,1)
\]
for any $\gamma,\lambda_\ell >0$. Thus, Theorem~\ref{the.gam.zer} yields that
\[
\lambda_\ell\;\, C^*(\gamma,\lambda_\ell)\stackrel{\rm d}{\to } Z\qquad\mbox{as}\quad \lambda_\ell\to \infty\,,
\]
where $Z\sim \mbox{Exp}(2)$.

\subsection{Asymptotic Weibull distribution for $\kappa\to \infty$}

In this section we assume that $T$ is a stationary and isotropic
random tessellation which is mixing. Furthermore, we assume that
 \begin{equation}\label{int.con.deb}
\E  \,\nu_1^2(\partial \Xi^*) <\infty\,,
\end{equation}
where $\nu_1(\partial \Xi^*)$ denotes the circumference of the typical cell $\Xi^*$ of $T$.

We investigate the asymptotic
behavior of the distribution of $C^*=C^*(\gamma,\lambda_\ell)$ for
$\kappa\to \infty$, where $\gamma\to \infty$ and $\lambda_\ell\to
0$ such that $\lambda_\ell\gamma=\lambda$ is fixed. In particular,
we show that $C^*$ converges in distribution to $\xi Z$, where $\xi \ge 1$ is a certain constant which is multiplied
by the (random) Euclidean distance $Z$ from the origin to the nearest point of a stationary
Poisson process of intensity $\lambda$.  Then, it is easy to see
that $Z$ as well as $\xi Z$ have Weibull distributions.
\begin{thm}\label{the.gam.inf}
    Let
    $Z\sim {\rm Wei}(\lambda\pi,2)$ for some $\lambda>0$. Then there exists a constant $\xi\ge 1$ such that
    \begin{equation}\label{con.ceh.iks}
        C^*(\gamma,\lambda_\ell)\stackrel{\rm d}{\to } \xi Z\qquad\mbox{as}\quad \kappa\to \infty\
    \end{equation}
    provided that $\gamma\to \infty$ and $\lambda_\ell\to 0$ with $\lambda_\ell\gamma=\lambda$,
 where $\xi Z\sim Wei(\lambda\pi/\xi^2,2)$.
\end{thm}

The {\em proof} of Theorem~\ref{the.gam.inf} is split into several
steps. We first show in Lemma~\ref{lem.cox.conv} that under the
Palm probability measure $\mathbb{P}^*_{X_L}$, the Euclidean
distance $|\widetilde{X}_{H,0}|$ from the origin to the nearest
point $\widetilde{X}_{H,0}$  of the point process
$\widetilde{X}_H=\{\widetilde{X}_{H,n}\}$ of HLC converges in
distribution to the corresponding characteristic of a stationary
Poisson process with intensity $\lambda$. Furthermore, in
Lemma~\ref{lem.conv.prob}, we show that for some constant $\xi \ge 1$
the difference between $\xi|\widetilde{X}_{H,0}|$ and
the shortest path length $C^*=C^*(\gamma,\lambda_\ell)$  from the
origin to $\widetilde{X}_{H,0}$ along the edge set
$\widetilde{T}^{(1)}_\gamma$ converges in probability to zero.
Then, combining the results of Lemmas~\ref{lem.cox.conv} and
\ref{lem.conv.prob}, the assertion of Theorem~\ref{the.gam.inf}
follows.

\section{Proof of Theorem~\ref{the.gam.inf}}\label{sec.fou}

\subsection{Some auxiliary results on convergence of point processes}

In the proofs of Lemmas~\ref{lem.alf.ren} and \ref{lem.cox.conv}
 which will be given below, we use two
classic results regarding the convergence in distribution of point
processes, see e.g. \cite{Daley0307,Kal86,KMM78}. Note that a
sequence of point processes $X^{(1)},X^{(2)},\ldots$ in $\R^2$ is
said to converge in distribution to a point process $X$ in $\R^2$
if
\[
\lim_{m\to\infty}\mathbb{P}(X^{(m)}(B_1)=i_1,\ldots,X^{(m)}(B_k)=i_k)=\mathbb{P}(X(B_1)=i_1,\ldots,X(B_k)=i_k)
\]
for any $k\ge 1$, $i_1,\ldots,i_k\ge 0$ and
 for all finite sequences of bounded sets $B_1,\ldots,B_k\in\mathcal{B}^2$ which satisfy the condition
$\mathbb{P}(X(\partial B_j)>0)=0$ for each $j=1,\ldots,k$. In this case we shortly write $X^{(m)}\Longrightarrow X$.

Let $X=\{X_n\}$ be an arbitrary ergodic point process
in $\R^2$ with $\mathbb{P}(X(\R^2)=0)=0$, and let $\lambda\in(0,\infty)$ denote the intensity of $X$.
Then, the following limit theorem for independently thinned and appropriately
re-scaled versions of $X$ is true.
For each $c\in(0,1)$, let $X^{(c)}$ denote a point process which arises from $X$ by independent thinning, where each atom $X_n$ of $X$
is deleted with probability $1-c$ (and ,,survives'' with probability $c$). Furthermore, let $Y^{(c)}$ be a re-scaled version of $X^{(c)}$,
where $Y^{(c)}(B)=X^{(c)}(B/\sqrt{c})$ for each $B\in\mathcal{B}^2$. Then, for each $c\in(0,1)$, the point process
$Y^{(c)}$ is  stationary with the same intensity $\lambda$ as $X$, and
\begin{equation}\label{the.alf.ren}
Y^{(c)}\Longrightarrow Y\qquad\mbox{if $c\to 0$,}
\end{equation}
where $Y$ is a stationary Poisson process in $\R^2$ with intensity
$\lambda$, see e.g. Section~11.3 of \cite{Daley0307} or
Theorem~7.3.1 in \cite{KMM78}. Moreover, the following continuity
property of Palm distributions holds. Let
$X,X^{(1)},X^{(2)},\ldots$ be stationary point processes in $\R^2$
such that $\mathbb{P}(X(\R^2)=0)= \mathbb{P}(X^{(m)}(\R^2)=0)=0$
for each $m\ge 1$ and let
 $\lambda,\lambda^{(1)},\lambda^{(2)},\ldots$ denote the intensity of $X,X^{(1)},X^{(2)},\ldots$, respectively.
If $\lambda^{(m)}=\lambda$ for each $m\ge 1$ and
$X^{(m)}\Longrightarrow X$ as $m\to\infty$, then
\begin{equation}\label{the.con.pal}
Y^{(m)}\Longrightarrow Y\qquad\mbox{as $m\to \infty$,}
\end{equation}
where $Y,Y^{(1)},Y^{(2)},\ldots$ are point processes in $\R^2$ whose distribution is equal to the
Palm distribution of $X,X^{(1)},X^{(2)},\ldots$, respectively, see e.g. Proposition~10.3.6   in \cite{KMM78}.

\subsection{Euclidean distance from the typical LLC to its closest HLC}\label{sub.euc.dis}

Throughout this section we assume that the underlying tessellation $T$ is ergodic.
In order to prove that the Euclidean distance $|\widetilde{X}_{H,0}|$ from the typical LLC to its closest HLC is asymptotically Weibull distributed, we first show
that the (stationary) Cox process $X_H$ converges in distribution to a homogeneous Poisson process if $\kappa\to \infty$
provided that $\lambda_\ell\gamma=\lambda$ is constant.

\begin{lem}\label{lem.alf.ren}
If $\kappa\to \infty$, where
 $\lambda_\ell\gamma=\lambda$ for some constant $\lambda\in(0,\infty)$,
then $X_H\Longrightarrow Y$, where $Y$ is a stationary Poisson process with intensity $\lambda$.
\end{lem}
\begin{proof}
For each $\gamma>1$, let $X_H=X_H(\gamma)$ denote the Cox process of HLC with parameters $\gamma$ and $\lambda_\ell$,
where $\lambda_\ell=\lambda/\gamma$ for some constant $\lambda\in(0,\infty)$. Note that the Cox process $X_H(\gamma)$ can be obtained from $X_H(1)$ by independent thinning with
survival probability $c=1/\gamma$ and by subsequent re-scaling with scaling factor $\sqrt{1/\gamma}$. Furthermore, the Cox process
$X_H(1)$ is ergodic, since $T$ is ergodic. Thus, using (\ref{the.alf.ren}), we get that
$X_H(\gamma)\Longrightarrow Y$ as $\gamma\to\infty$.
\end{proof}

\begin{lem}\label{lem.cox.conv}
    Let
    $Z\sim {\rm Wei}(\lambda\pi,2)$ for some $\lambda>0$. Then
     $|\widetilde{X}_{H,0}|\stackrel{\rm d}{\to } Z$ as $\quad \kappa\to
     \infty$
    provided that $\gamma\to \infty$ and $\lambda_\ell\to 0$ such that $\lambda_\ell\gamma=\lambda$.
\end{lem}
\begin{proof}
Let $X_H^*(\gamma)$ be a point process in $\R^2$ whose distribution is equal to the Palm distribution of $X_H=X_H(\gamma)$.
Furthermore, let $Y$ be a stationary Poisson process with intensity $\lambda$. Note that  the distribution of  $Y+\delta_o$
is then equal to the Palm distribution of $Y$, see e.g. Proposition 13.1.VII in \cite{Daley0307}.
Thus, using (\ref{the.con.pal}), Lemma~\ref{lem.alf.ren} gives that
\begin{equation}\label{con.pal.ver}
X_H^*(\gamma)\Longrightarrow Y+\delta_o
\end{equation}
as $\gamma\to\infty$ and $\lambda_\ell\to 0$, where $\lambda_\ell\gamma=\lambda$.
Since $X_L$ and $X_H$ are Cox processes concentrated on $T_\gamma^{(1)}$
    which are conditionally independent given $T_\gamma^{(1)}$, we get that $\widetilde{X}_H+\delta_0$ and the Palm version $X_H^*$
    of $X_H$ have the same distributions. This is an easy consequence of the representation formula for the Palm distribution
of stationary Cox processes, see e.g. Section~5.2 in \cite{Stoyan95}. In particular, this gives that for each $r>0$
\begin {eqnarray*}
\lim_{\gamma\to\infty}  \mathbb{P}(|\widetilde{X}_{H,0}|>r) &=& \lim_{\gamma\to\infty}
 \mathbb{P}(\widetilde{X}_H(B(o,r))=0)\\
&=&  \lim_{\gamma\to\infty} \mathbb{P}((\widetilde{X}_H+\delta_o)(B(o,r))=1)\\
&=&  \lim_{\gamma\to\infty}  \mathbb{P}(X_H^*(B(o,r))=1)\\
&=&  \mathbb{P}((Y+\delta_0)(B(o,r))=1)\\
&=&  \mathbb{P}(Y(B(o,r))=0) \,,
    \end{eqnarray*}
where we used (\ref{con.pal.ver}) in the last but one equality.
Thus, for each $r>0$,
\[
\lim_{\gamma\to\infty}  \mathbb{P}(|\widetilde{X}_{H,0}|>r)\;=\;  \mathbb{P}(Y(B(o,r))=0)
\;=\; \exp(-\lambda\pi r^2)\,,
\]
which means that $|\widetilde{X}_{H,0}|\stackrel{\rm d}{\to } Z\sim {\rm Wei}(\lambda\pi,2)$.
\end{proof}

\subsection{Shortest path length vs. scaled Euclidean distance}

In this section we assume that $T$ is a stationary and isotropic
random tessellation which is mixing. Furthermore, we assume that the integrability condition (\ref{int.con.deb}) is
satisfied.
Then, we can show that for
some constant $\xi \ge 1$ the difference between
$\xi|\widetilde{X}_{H,0}|$ and the shortest path length
$C^*=C^*(\gamma,\lambda_\ell)$  from the origin to
$\widetilde{X}_{H,0}$ along the edge set
$\widetilde{T}^{(1)}_\gamma$ converges in probability to zero. In
order to show this we need the following auxiliary result.

\begin{lem}\label{lem.exp.conv}
Let $\widetilde{T}_{\gamma, \xi, \varepsilon}^{(1)}=\big\{u\in
\widetilde{T}_\gamma^{(1)}: \big|c(u)-\xi|u|\big| <
\varepsilon\big\}$, where $\xi\ge 1$ is some constant and $c(u)$ denotes the length of the
shortest path from $u$ to the origin along the edges of
$\widetilde{T}_\gamma^{(1)}$. If $\gamma\to \infty$ and
$\lambda_\ell\to 0$, where $\lambda_\ell\gamma=\lambda$ is fixed,
then there exists $\xi\ge 1$ such that for each
$\varepsilon>0$ and $r>0$
    \begin{equation}\label{exp.con.aux}
        \lim_{\gamma\rightarrow \infty}\E
        \exp\Big(-\frac{\lambda}{\gamma}\,
\nu_1\big(\widetilde{T}^{(1)}_\gamma\backslash\widetilde{T}^{(1)}_{\gamma, \xi, \varepsilon}\cap
B(o,r)\big)\Big)=1\,.
    \end{equation}
\end{lem}

The {\em proof} of this lemma is postponed to the Appendix. Now,
using Lemma~\ref{lem.exp.conv}, we are able to complete the proof
of Theorem~\ref{the.gam.inf} by showing that the following is
true.

\begin{lem}\label{lem.conv.prob}
If $\gamma\to \infty$ and $\lambda_\ell\to 0$ such that
$\lambda_\ell\gamma=\lambda$, then there is a constant $\xi\ge 1$ with
    $ C^*(\gamma,\lambda_\ell)-\xi |\widetilde{X}_{H,0}|
        \stackrel{\rm P}{\to } 0$, where $ \stackrel{\rm P}{\to }$ denotes convergence in
    probability.
\end{lem}

\begin{proof} We have to show that there exists a constant $\xi\ge 1$ such that
for any $\varepsilon>0$ and $\delta > 0$ we can choose
$\gamma_0>0$ with
\begin{eqnarray*}
    \mathbb{P}\big(\big|C^*-\xi |\widetilde{X}_{H,0}|\big|>\varepsilon\big)\le\delta
\end{eqnarray*}
for all $\gamma>\gamma_0$. Note that
\begin{eqnarray*}
    \lefteqn{\mathbb{P}\big(\big|C^*-\xi |\widetilde{X}_{H,0}|\big|>\varepsilon\big)}\\
    &&=\mathbb{P}\big(\big|C^*\!\!- \xi |\widetilde{X}_{H,0}|\big|>\varepsilon,|\widetilde{X}_{H,0}| \le r\big)
    +\mathbb{P}\big(\big|C^*\!\!-\xi |\widetilde{X}_{H,0}|\big|>\varepsilon,|\widetilde{X}_{H,0}| > r\big)\,,
\end{eqnarray*}
where $r>0$ is an arbitrary fixed
number. Since 
\[
\mathbb{P}\big(|\widetilde{X}_{H,0}|
>r\big)\longrightarrow e^{-\lambda\pi r^2}\quad\mbox{as}\quad \gamma \longrightarrow \infty\,,
\]
see Lemma~\ref{lem.cox.conv}, we can choose $r>0$ such that
$\mathbb{P}\big(|\widetilde{X}_{H,0}| > r\big)< \delta/2$ for all
$\gamma
> 0$ sufficiently large. Thus, it is enough to show that there
exists $\gamma_0>0$ such that $ \mathbb{P}\big(\big|C^*\!\!-
\xi |\widetilde{X}_{H,0}|\big|>\varepsilon,|\widetilde{X}_{H,0}| \le r\big)
\le\delta/2$ for all $\gamma>\gamma_0$. Let
$\widetilde{N}=\widetilde{X}_H(B(o,r))$ denote the number of
points of $\widetilde{X}_H$ in $B(o,r)$. Then we have
\begin{eqnarray*}
    \lefteqn{ \mathbb{P}\big(\big|C^*\!\!-\xi
|\widetilde{X}_{H,0}|\big|>\varepsilon,|\widetilde{X}_{H,0}| \le r\big)}\\
    &&\le \E \left(\sum\limits_{k=1}^\infty
\mathbb{P}(\widetilde{N}=k\mid\widetilde{T}_\gamma)\;\mathbb{P}\Bigl(\max\limits_{i=1,\dots,k}\big(\big|c(Y_i)
-\xi |Y_i|\big|\big)>\varepsilon  \,\,\Big|\,\, \widetilde{T}_\gamma, \widetilde{N}=k\Bigr)\right)\\
    &&=\E \left(\sum\limits_{k=1}^\infty \mathbb{P}(\widetilde{N}=k\mid\widetilde{T}_\gamma)\,\Bigl(1-\mathbb{P}\big(\big|c(Y_1)
    -\xi|Y_1|\big|\le\varepsilon \,\,\big|\,\, \widetilde{T}_\gamma\bigr)^k\Bigr)\right)\,,
\end{eqnarray*}
where the points $Y_1,\dots,Y_k$ are conditionally independent and
identically distributed according to $\nu_1\big(\,\cdot\,\cap\,
\widetilde{T}_\gamma^{(1)}\cap
B(o,r)\big)/\nu_1\big(\widetilde{T}_\gamma^{(1)}\cap B(o,r)\big)$ for given
$\widetilde{T}_\gamma$ and $\widetilde{N}=k$. In particular, for
the conditional probability in the latter expression, we have
\begin{eqnarray*}
    \mathbb{P}\big(\big|c(Y_1)
    -\xi|Y_1|\big|\le\varepsilon \mid \widetilde{T}_\gamma\big)
    &=&\frac{\displaystyle\int_{\widetilde{T}_\gamma^{(1)}\cap B(o,r)}\ind_{[-\varepsilon,\varepsilon]}(c(u)-\xi|
    u|)\,\nu_1(du)}{\nu_1(\widetilde{T}_\gamma^{(1)}\cap
B(o,r))}
    \\
    &=&\frac{\nu_1\big(\widetilde{T}_{\gamma, \xi, \varepsilon}^{(1)}\cap B(o,r)\big)}{\nu_1\big(\widetilde{T}_\gamma^{(1)}\cap
    B(o,r)\big)}\;.
\end{eqnarray*}
Using that $\widetilde{N}\sim Poi(\widetilde{\lambda})$ with $\widetilde{\lambda}=\lambda_\ell\nu_1\big(\widetilde{T}_\gamma^{(1)}\cap B(o,r)\big)$ given
$\widetilde{T}_\gamma$, we get
\begin{eqnarray*}
    \lefteqn{\sum\limits_{k=1}^\infty \mathbb{P}(\widetilde{N}=k\mid\widetilde{T}_\gamma)\,\Bigl(1-\mathbb{P}\big(\big|c(Y_1)
    -\xi|Y_1|\big|\le\varepsilon \,\,\big|\,\, \widetilde{T}_\gamma\bigr)^k\Bigr)}\hspace{5cm}\\&&=
    1-\sum\limits_{k=0}^\infty e^{-\widetilde{\lambda}}
    \frac{\widetilde{\lambda}^k}{k!}
    \Bigl( \frac{\lambda_\ell\nu_1\big(\widetilde{T}_{\gamma, \xi, \varepsilon}^{(1)}\cap B(o,r)\big)}{\widetilde{\lambda}}\Bigr)^k\\
    &&=1- \sum\limits_{k=0}^\infty e^{-\widetilde{\lambda}}
    \frac{1}{k!}
    \Bigl( \lambda_\ell\nu_1\big(\widetilde{T}_{\gamma, \xi, \varepsilon}^{(1)}\cap B(o,r)\big)\Bigr)^k\\
    &&=1-e^{-\lambda_\ell\bigr( \nu_1(\widetilde{T}_\gamma^{(1)}\cap B(o,r)) - \nu_1(\widetilde{T}_{\gamma, \xi, \varepsilon}^{(1)}\cap B(o,r))\bigl)}\,.
\end{eqnarray*}
Thus we have
\[
    \lim_{\gamma\to \infty}\mathbb{P}\big(\big|C^*\!-\xi
|\widetilde{X}_{H,0}|\big|>\varepsilon,|\widetilde{X}_{H,0}| \le
r\big)\,\le\,1- \lim_{\gamma\rightarrow \infty}\E
        \exp\Big(-\frac{\lambda}{\gamma}\,
\nu_1\big(\widetilde{T}^{(1)}_\gamma\backslash\widetilde{T}^{(1)}_{\gamma, \xi, \varepsilon}\cap
B(o,r)\big)\Big)
 \,.
\]
Using Lemma~\ref{lem.exp.conv} this gives that $ \lim_{\gamma\to
\infty} \mathbb{P}\big(\big|C^*-\xi
|\widetilde{X}_{H,0}|\big|>\varepsilon,|\widetilde{X}_{H,0}| \le
r\big)\;=\;0$, which completes the proof.
\end{proof}


\section{Examples}\label{sec.sec.fiv}

Recall that in Theorem~\ref{the.gam.inf} we assumed that the
underlying tessellation $T$ is stationary and isotropic. The
examples of tessellations discussed in the present section
obviously possess these properties. Furthermore, we assumed in
Theorem~\ref{the.gam.inf} that $T$ is mixing and fulfills the
integrability condition (\ref{int.con.deb}). We first show that
the mixing condition is satisfied for a wide class of
tessellations. Moreover, we also show that
(\ref{int.con.deb}) is true for these tessellations.

The tessellation models considered
in the literature focus mainly on PLT and PVT as well as on
Poisson-Delaunay tessellations (PDT), on iterated tessellations
constructed from these basic tessellations of Poisson type and on STIT tessellations, see
e.g. \cite{Ald08}--\cite{Bac96}, \cite{FL08}--\cite{Glo09}, \cite{Nagel05} and
\cite{Voss08b}--\cite{Weiss99}. Here, we assume that an iterated
tessellation is either a $T_I/T_{II}$-superposition or a
$T_I/T_{II}$-nesting of tessellations $T_I$ and $T_{II}$ as
defined e.g. in \cite{Ba02,Mai03,Weiss99}.
Note that the edge set of a $T_I/T_{II}$-superposition is given
by the union $T_I^{(1)}\cup T_{II}^{(1)}$, where $T_I$ and $T_{II}$ are independent.
Furthermore, a $T_I/T_{II}$ nesting is constructed by subdividing each cell
of $T_I$ by independent copies of $T_{II}$.
We show that for these
important models Theorem~\ref{the.gam.inf} can be applied.
Furthermore, if $T$ is a PLT or a
$T_I/T_{II}$-superposition/nesting with $T_I$ being a PLT, then we
can even calculate the constant $\xi$ explicitly that appears in
Theorem~\ref{the.gam.inf}. On the other hand, if $T$ is a PDT,
 we get an upper bound for $\xi$.

\subsection{Mixing tessellations}\label{sub.mix.tes}

In order to apply Theorem~\ref{the.gam.inf} we have to show that
the underlying tessellation $T$ is mixing, where we will use the
following criterion to show that a stationary random closed set is
mixing.
\begin{lem}\label{lem.fiv.one}
    A stationary random closed set $\Xi$ in $\R^2$ is mixing if and only if
    \begin{equation}\label{eq.mix}
        \lim_{|x|\rightarrow \infty}\mathbb{P}(\Xi \cap C_1=\emptyset,\, \Xi\cap (C_2+x)=\emptyset)=
        \mathbb{P}(\Xi \cap C_1=\emptyset)\;\mathbb{P}(\Xi \cap C_2=\emptyset)
    \end{equation}
    for all $C_1, C_2\in\mathcal{R}$, where $\mathcal{R}$ is the family of all subsets of $\R^2$
    which are finite unions of closed balls with rational radii and centres
    with rational coordinates.
\end{lem}

 Note that the statement of Lemma~\ref{lem.fiv.one} is essentially
 Lemma~4 in \cite{Hei92}, see also Theorem~9.3.2 in \cite{ScW08},
    where the (stronger) condition is considered that (\ref{eq.mix})  holds for all
    compact sets $C_1, C_2\subset\R^2$. However, it is easy to see that it suffices to assume that
    (\ref{eq.mix}) holds for the separating class $\mathcal{R}$; see also Section~1.4 of \cite{Mol05}.
    To make this clear, we only have to show that
     $   \mathcal{E}=\{\mathcal{F}^{C_0}_{C_1,\dots,C_k}:C_0,\dots,C_k\in \mathcal{R}^\prime,k\ge
     0\}$
    is a semi-algebra which generates $\mathcal{B}(\mathcal{F})$,
    where $\mathcal{R}^\prime=\mathcal{R}\cup\emptyset$ and
    \[
    \mathcal{F}^{C_0}_{C_1,\dots,C_k}=\{F\in\mathcal{F}:
    F\cap C_0=\emptyset,F\cap C_1\not=\emptyset,\dots,F\cap C_k\not=\emptyset\}\,.
    \]
    Note that the family $\mathcal{R}^\prime$ is union-stable. Thus, by Lemma~2.2.2 in
    \cite{ScW08}, we get
     that $\mathcal{E}$ is a semi-algebra. Moreover, let $G\subset\R^2$ denote an open set, then $G=\bigcup_{i=1}^\infty C_i$
    for some $C_1,C_2,\ldots\in\mathcal{R}^\prime$ and $\mathcal{F}_G=\{F\in\mathcal{F}:F\cap G\not=\emptyset\}=
    \bigcup_{n=1}^\infty \mathcal{F}_{\bigcup_{i=1}^n C_i}$, thus $\mathcal{F}_G\in\sigma(\mathcal{E})$.
    Since $\{\mathcal{F}_G:G\subset\R^2\,\,\mbox{open}\}$
    generates $\mathcal{B}(\mathcal{F})$, we get that $\sigma(\mathcal{E})=\mathcal{B}(\mathcal{F})$.
    Now the statement of Lemma~\ref{lem.fiv.one} can be proven by exactly the same arguments used in the proof of
	Lemma~4 in \cite{Hei92}.

It is well known that $T$ is mixing  if $T$ is a PDT, PVT and PLT,
respectively, see e.g. Chapter~10.5 in \cite{ScW08}, and recently it was shown
that STIT tessellations are mixing (\cite{Lac09}). Furthermore,
using Lemma~\ref{lem.fiv.one}, we can show that $T$ is  mixing if
$T$ is an iterated tessellation constructed from these basic
tessellations of Poisson type.

\begin{lem}\label{lem.mix.ite}
The tessellation $T$ is mixing  if $T$ is a
$T_I/T_{II}$-superposition of two mixing  tessellations $T_I$ and
$T_{II}$, or a $T_I/T_{II}$-nesting of a mixing initial
tessellation $T_I$ and any stationary component
    tessellation $T_{II}$.
\end{lem}

\begin{proof}
Suppose first that $T$ is a $T_I/T_{II}$-superposition. Then, for
any $C_1, C_2\in\mathcal{R}$
\begin{eqnarray*}
\lefteqn{\mathbb{P}(T^{(1)} \cap C_1=\emptyset,\,  T^{(1)}\cap
(C_2+x)=\emptyset)}\\
&=& \mathbb{P}(T_I^{(1)} \cap C_1=\emptyset,\, T_I^{(1)}\cap
(C_2+x)=\emptyset,\,T_{II}^{(1)} \cap C_1=\emptyset,\,
T_{II}^{(1)}\cap (C_2+x)=\emptyset)\\
&=& \mathbb{P}(T_I^{(1)} \cap C_1=\emptyset,\, T_I^{(1)}\cap
(C_2+x)=\emptyset)\;  \mathbb{P}(T_{II}^{(1)} \cap
C_1=\emptyset,\, T_{II}^{(1)}\cap (C_2+x)=\emptyset)\,,
\end{eqnarray*}
since $T_I$ and $T_{II}$ are independent. Thus, using
Lemma~\ref{lem.fiv.one}, we get that $T$ is mixing if $T_I$ and
$T_{II}$ are mixing. Let now $T$ be a $T_I/T_{II}$-nesting and
assume that $C_1=\cup_{j=1}^n B_j,C_2=\cup_{j=n+1}^{n+m} B_j$ for
closed balls $B_1,\dots,B_{n+m}\subset\R^2$ with rational radii
and centres with rational coordinates. Let $\Xi_1,\Xi_2,\ldots$ be
the cells of the initial tessellation $T_I=\{\Xi_n\}$,  let $D$
denote the family of all decompositions of the index set
$\{1,\dots,n+m\}$ into nonempty subsets, and for
$J=\{J_1,\dots,J_k\}\in D$ consider the set
\begin{equation}\label{def.aah.jot}
 A_J(x)=\{\cup_{j\in
J_i}(B_j+x\ind_{\{j>n\}})\subset \inter\Xi_{j_i},\; i=1,\dots,k,\;
\Xi_{j_i}\not=\Xi_{j_l}\,\mbox{for}\, j_i\not=j_l\}\,,
\end{equation}
 i.e., each
of the sets $\cup_{j\in J_i}(B_j+x\ind_{\{j>n\}})$ is contained in
a different cell of $T_I$. Using this notation we get
\begin{eqnarray*}
  \lefteqn{\lim_{|x|\rightarrow \infty}\mathbb{P}(T^{(1)} \cap C_1=\emptyset,\, T^{(1)}\cap
  (C_2+x)=\emptyset)}\\
    &=&\sum_{J\in D}\lim_{|x|\rightarrow \infty}\mathbb{P}(T^{(1)} \cap C_1=\emptyset,\, T^{(1)}\cap
  (C_2+x)=\emptyset,\, A_J(x))\,.
\end{eqnarray*}
Since the cells $\Xi_1,\Xi_2,\ldots$ of $T_I$ are finite with
probability $1$, we have
\begin{eqnarray*}
    \lim_{|x|\rightarrow \infty} \mathbb{P}(T^{(1)} \cap C_1=\emptyset,\, T^{(1)}\cap
  (C_2+x)=\emptyset,\, A_J(x)) &=& 0
\end{eqnarray*}
if there are $i\le n$ and $j>n$ with $i,j\in J_l\in J$.  On the other hand, suppose that
$J=\{J_1,\dots,J_k\}$ is a decomposition of $\{1\dots,n+m\}$ with
$J_i\subset\{1,\dots,n\}$ for $i=1,\dots,l$ and
$J_i\subset\{n+1,\dots,n+m\}$ for $i=l+1,\dots,k$. Then we get
that
\begin{eqnarray*}
    \lefteqn{\mathbb{P}(T^{(1)} \cap C_1=\emptyset,\, T^{(1)}\cap
  (C_2+x)=\emptyset,\, A_J(x))}\\
    &&=\mathbb{P}(A_J(x),\, B_{J_i}\cap
    T^{(1)}_{II,i}=\emptyset,\,i=1,\dots,l,\,
    B_{J_i}+x\cap T^{(1)}_{II,i}=\emptyset,\,i=l+1,\dots,k)\,,
\end{eqnarray*}
where $B_{J_i}=\cup_{j\in J_i}B_j$ and $T_{II,1},\dots,T_{II,k}$ are independent copies of $T_{II}$
which are independent of $T_I$.
Thus we have
\begin{eqnarray*}\label{iterated.mix.2}\nonumber
    \lefteqn{\mathbb{P}(A_J(x),\, B_{J_i}\cap
    T^{(1)}_{II,i}=\emptyset,\,i=1,\dots,l,\,
    B_{J_i}+x\cap T^{(1)}_{II,i}=\emptyset,\,i=l+1,\dots,k)}\\
    &&=\mathbb{P}(A_J(x))\;\mathbb{P}(B_{J_i}\cap T^{(1)}_{II,i}=\emptyset,\,i=1,\dots,l)\;
    \mathbb{P}(B_{J_i}\cap T^{(1)}_{II,i}=\emptyset,\,i=l+1,\dots,k)
    \,.
\end{eqnarray*}
Moreover, since $T_I$ is mixing, we get
\begin{eqnarray*}\label{iterated.mix.3}
    \lim_{|x|\rightarrow \infty}\mathbb{P}(A_J(x))&=&\mathbb{P}(A_{J'}(o))\;\mathbb{P}(A_{J''}(o))\,,
\end{eqnarray*}
where $J'=\{J_1,\dots,J_l\}$ and $J''=\{J_{l+1},\dots,J_k\}$ are
the decompositions of $\{1,\dots,n\}$ and $\{n+1,\dots,n+m\}$,
respectively, induced by $J$, and $A_{J'}(o)$ resp. $ A_{J''}(o)$
are defined analogously to (\ref{def.aah.jot}). Summarizing the
above considerations, we get
\begin{eqnarray*}
    \lefteqn{\lim_{|x|\rightarrow \infty}\mathbb{P}(T^{(1)} \cap C_1=\emptyset,\, T^{(1)}\cap
  (C_2+x)=\emptyset,\, A_J(x)) }\\
    &=&\mathbb{P}(A_{J'}(o),\,B_{J_i}\cap T^{(1)}_{II,i}=\emptyset,\,i=1,\dots,l)\\
    && \hspace{5cm}
   \times\;\; \mathbb{P}(A_{J''}(o),\,B_{J_i}\cap T^{(1)}_{II,i}=\emptyset,\,i=l+1,\dots,k)\\
    &=&\mathbb{P}(T^{(1)}\cap C_1=\emptyset,A_{J'}(o))\;\mathbb{P}(T^{(1)}\cap C_2=\emptyset,A_{J''}(o))\,,
\end{eqnarray*}
which yields
\begin{eqnarray*}
\lefteqn{\lim_{|x|\rightarrow \infty}\mathbb{P}(T^{(1)} \cap
C_1=\emptyset,\, T^{(1)}\cap
  (C_2+x)=\emptyset)}\\
    &=&\sum_{J\in D}\lim_{|x|\rightarrow \infty}\mathbb{P}(T^{(1)} \cap C_1=\emptyset,\, T^{(1)}\cap
  (C_2+x)=\emptyset,\, A_J(x))\\
    &=&\sum_{J'\in D'}\sum_{J''\in D''}\mathbb{P}(T^{(1)}\cap C_1=\emptyset,A_{J'}(o))\;\mathbb{P}(T^{(1)}\cap C_2=\emptyset,A_{J''}(o))\\
    &=&\mathbb{P}(T^{(1)}\cap C_1=\emptyset)\;\mathbb{P}(T^{(1)}\cap C_2=\emptyset)\,,
\end{eqnarray*}
where $D'$, $D''$ is the family of all decompositions of
$\{1,\dots,n\}$ and $\{n+1,\dots,n+m\}$, respectively. Thus, by
Lemma~\ref{lem.fiv.one}, the nested tessellation $T$ is mixing.
\end{proof}

\subsection{Integrability condition (\ref{int.con.deb})}\label{sub.int.con}

The next result provides several classes of stationary
tessellations such that the second moment of the circumference of
their typical cell is finite, where $R(\Xi)$ denotes the radius of
the minimal ball containing the random convex polygon $\Xi$.

\begin{lem}\label{lem.exa.fin}
    If $T$ is a PVT, PDT, PLT and STIT tessellation, respectively, then $\E R^2(\Xi^*)  <
    \infty$ and, consequently,
\begin{equation}\label{fin.sec.mom}
\E \nu_1^2(\partial \Xi^*)  < \infty\,.
\end{equation}
Moreover, $(\ref{fin.sec.mom})$ holds
    if $T$ is a a
$T_I/T_{II}$-superposition/nesting such that
\begin{equation}\label{fin.ite.mom}
    \max \{\E R^2( \Xi_I^*), \E R^2(\Xi_{II}^*)\}<\infty\,,
\end{equation}
     where $\Xi^*_I$
    and $\Xi^*_{II}$ is the typical cell of $T_I$ and $T_{II}$,
    respectively.
\end{lem}
\begin{proof}
Note that
\begin{equation}\label{ine.mu.er}
 \E\nu_1^2(\partial \Xi^*) \le 4\pi^2 \E R^2(\Xi^*)
\end{equation}
  holds for the typical cell $\Xi^*$ of
any stationary tessellation $T$. Furthermore, if $T$ is a PDT,
then it is well known that $\E R^2(\Xi^*)<\infty$. This result
goes back to \cite{Mil74}, see also Theorem~7.5 in \cite{Moe89}
and Theorem~10.4.4 in \cite{ScW08}. Similarly, it is well known
that $\E R^2(\Xi^*)<\infty$ holds if $T$ is a PVT or PLT, see e.g.
\cite{Cal02}. Since the interior of the typical cell of a STIT tessellation and a PLT
have the same distribution (\cite{Nagel05}), it is clear that $\E R^2(\Xi^*)<\infty$ also for STIT tessellations.
    If $T=T_I/T_{II}$ is an iterated tessellation with cell intensity $\lambda_T$,
    then we can use Proposition~3.1 in \cite{Mai04}
    and Campbell's theorem in order to get
    \begin{eqnarray*}
        \E \nu^2_1(\partial \Xi^*)
        &=& \frac{\lambda_I}{\lambda_T}\; \E\Bigl(\sum_{\Xi_i\in T_{II}}
            \nu^2_1(\partial (\Xi_i\cap \Xi_I^*))\,\ind_{\{\Xi_i\cap \Xi_I^*\not=\emptyset\}}\Bigr)\\
        &=&\frac{\lambda_I\lambda_{II}}{\lambda_T}\; \E\int_{\R^2}
            \nu^2_1(\partial (\Xi_{II}^*+x\cap \Xi_I^*))\,\ind_{\{\Xi_{II}^*+x\cap
            \Xi_I^*\not=\emptyset\}}\,\nu_2(dx)\,,
\end{eqnarray*}
where $\lambda_I, \lambda_{II}$ and $\Xi^*_I, \Xi^*_{II}$ denote the cell intensities and the typical cells, respectively,
of $T_I$ and $T_{II}$. Note that we can assume that $\Xi^*_I$ and $\Xi^*_{II}$ are independent random convex bodies.
Since $\nu_1^2(\partial (\Xi^*_{II}+x\cap \Xi^*_I))\le \min\{\nu_1^2(\partial\Xi^*_{I}), \nu_1^2(\partial\Xi^*_{II})\}$
we get
\begin{eqnarray*}
        \E \nu^2_1(\partial \Xi^*)
               &\le& \frac{\lambda_I\lambda_{II}}{\lambda_T}\; \E\Bigl(\min\{\nu^2_1(\partial \Xi_{I}^*),\nu^2_1(\partial \Xi_{II}^*)\}
            \nu_2(\check\Xi_{II}^*\oplus \Xi_I^*)\Bigr)\\
        &\le& \frac{4\pi\,\lambda_I\lambda_{II}}{\lambda_T}\; \E\Bigl(\min\{\nu^2_1(\partial \Xi_{I}^*),\nu^2_1(\partial \Xi_{II}^*)\}
            \;\max\{R^2(\Xi_{I}^*),R^2(\Xi_{II}^*)\}\Bigr)\,,
   \end{eqnarray*}
    where in the latter inequality we used that
    \[
        \nu_2(\check\Xi_{II}^*\oplus \Xi_I^*)\;\le\; \pi R^2(\check \Xi_I^*\oplus  \Xi_{II}^*)
        \;\le \; 4\pi\max\{R^2(\Xi_I^*),R^2(\Xi_{II}^*)\}\,.
    \]
Using (\ref{ine.mu.er}) and the independence of $\Xi^*_I$ and
$\Xi^*_{II}$, this gives
\begin{eqnarray*}
        \E \nu^2_1(\partial \Xi^*)
 &\le& \frac{4\pi^3\,\lambda_I\lambda_{II}}{\lambda_T}\; \E\Bigl(\min\{ R^2(\Xi_{I}^*),R^2(\Xi_{II}^*)\}
            \;\max\{R^2(\Xi_{I}^*),R^2(\Xi_{II}^*)\}\Bigr)\\
        &=& \frac{4\pi^3\,\lambda_I\lambda_{II}}{\lambda_T}\; \E R^2(\Xi_{I}^*)\,\E
        R^2(\Xi_{II}^*)\;<\;\infty\,,
  \end{eqnarray*}
 provided that (\ref{fin.ite.mom}) holds.
\end{proof}

\subsection{Asymptotic Weibull distribution of shortest path lengths}\label{sub.asy.sho}

In Sections~\ref{sub.mix.tes} and \ref{sub.int.con} we showed that the assumptions
of Theorem~\ref{the.gam.inf} are fulfilled  for
several classes of random tessellations $T$. Thus, we are now able
to apply Theorem~\ref{the.gam.inf} to these tessellations.

\begin{cor}\label{cor.tes.mod}
    Let $Z\sim Wei(\lambda\pi,2)$ and let $T$ be a
  PDT, PVT, PLT or STIT tessellation,
       or an  iterated tessellation $T=T_I/T_{II}$
        such that condition $(\ref{fin.ite.mom})$ is fulfilled,
        where $T$ is either
 \begin{enumerate}
        \item[\rm 1.]
          a superposition of two mixing   tessellations $T_I$ and
        $T_{II}$, or
 \item[\rm 2.]
        a nesting of a mixing initial tessellation $T_I$ and any stationary
component    tessellation $T_{II}$.
    \end{enumerate}
    Then
    $C^*\stackrel{{\rm d}}{\to} \xi Z$ for some constant $\xi\ge 1$ provided that 
    $\gamma\to \infty$ and $\lambda_\ell\to 0$ such that
    $\lambda_\ell\gamma=\lambda$.
    Furthermore, if $T$ is a PLT or a $T_I/T_{II}$-superposition/nesting, 
    where $T_I$ is a PLT, then $\xi=1$.
    If $T$ is a PDT, then $\xi \le 4/\pi\approx 1.27$.
\end{cor}

\begin{proof}
 The first part of the assertion follows from
Theorem~\ref{the.gam.inf} if the results of
Lemmas~\ref{lem.mix.ite} and \ref{lem.exa.fin} as well as the
comments immediately before Lemma~\ref{lem.mix.ite} are taken into
account.
Now we consider the cases that $T$ is a PLT, a $T_{I}/T_{II}$-superposition/nesting with a PLT $T_{I}$,
or a PDT. To begin with, let $T$ be a PLT with intensity 1.
Then, the edge set $\widetilde{T}_\gamma^{(1)}$ of the
tessellation $\widetilde{T}_\gamma$ introduced in
Section~\ref{subsec.ser.zon}  is generated by a random
sequence of lines $L_0,L_1,\dots$, where $L_1,L_2,\dots$ form the edge set $T_\gamma^{(1)}$ of the (stationary and isotropic)
PLT $T_\gamma$ and $L_0$ is an isotropic line through the origin $o$, which is independent of
$T_\gamma$. Thus we have
\begin{eqnarray*}
    \frac{1}{\gamma} \; \nu_1(\widetilde{T}_\gamma^{(1)}\backslash \widetilde{T}_{\gamma, \xi, \varepsilon}^{(1)}\cap B(o,r))
    &\le & \frac{1}{\gamma} \; \nu_1(T_\gamma^{(1)}\cap B(o,r))+\frac{2r}{\gamma}\;.
\end{eqnarray*}
Using Theorem~\ref{lem.uni.int}, together with Lemma~\ref{lem.dom.int},
this yields that the family of random variables $\{X_{\gamma,\xi},\,\gamma>0\}$ with
$X_{\gamma,\xi}=\nu_1(\widetilde{T}_\gamma^{(1)}\backslash
\widetilde{T}_{\gamma, \xi, \varepsilon}^{(1)}\cap B(o,r))/\gamma$ is
uniformly integrable since $\nu_1(T_\gamma^{(1)}\cap
B(o,r))/\gamma=\pi r^2 \nu_1(T^{(1)}\cap B(o,r\gamma))/\nu_2(B(o,r\gamma))$
converges to $r^2\pi$ in $L^1$ due to the fact that the PLT $T$ is mixing and, therefore, ergodic
(\cite{Daley0307}, Theorem 12.2.IV). Furthermore,
in Lemma~\ref{lem.exp.conv} we showed that there is a $\xi\ge 1$ such that the Laplace transform of
$X_{\gamma,\xi}$ converges to 1, which implies that
$X_{\gamma,\xi}\stackrel{\rm P}{\longrightarrow }0$ (\cite{Kal02},
Theorem~5.3). Thus, applying Theorem~\ref{lem.uni.int} again, we get that 
\begin{equation}\label{con.iks.zer}
\lim\limits_{\gamma\to\infty} \E \,X_{\gamma,\xi} = 0\,.
\end{equation}
However, if $T_\gamma^{(1)}$ gets denser, there are lines which intersect the line through $o$ close to $o$. Thus, all
points on these lines have approximately the direct connections as shortest paths which
can be used to show that $\xi=1$.  Suppose that $\xi > 1$ and let $r>2>\varepsilon >0$ with $\xi >
1+\varepsilon$. If the line $L_i$ intersects the  segment $L_{0,\varepsilon}$, where $L_{0,\varepsilon}=L_0\cap
B(o,\varepsilon/2)$, then for each $y\in L_i$ it holds that $0\le c(y)-|y| \le \varepsilon$ since the path from
$y$ to $o$ via the intersection point $L_i\cap L_{0,\varepsilon}$ is not longer than $|y|+\varepsilon$. Thus,
if $|y| > 2$ 
\[
    \big|c(y)-\xi|y|\big|\;=\; \big|c(y)-|y|-(\xi-1)|y|\big|
    \;\ge\; \varepsilon (|y|-1)\ge\varepsilon\,,
\]
which means that $y\in \widetilde{T}_\gamma^{(1)}\backslash \widetilde{T}_{\gamma, \xi, \varepsilon}^{(1)}$.
Furthermore, if $L_i\cap L_{0,\varepsilon}\not=\emptyset$, it is not difficult to see
that $\nu_1(L_i\cap B(o,r)\backslash B(o,2))\ge a$ for some constant $a>0$. These two observations lead to
\begin{eqnarray*}
   X_{\gamma,\xi}\;=\; \frac{1}{\gamma}  \nu_1(\widetilde{T}_\gamma^{(1)}\backslash \widetilde{T}_{\gamma, \xi, \varepsilon}^{(1)}\cap B(o,r))
    &\ge & \frac{1}{\gamma}\; \nu_1\Bigl(\bigcup_{i:L_i\cap L_{0,\varepsilon}\not=\emptyset }\{L_i\cap B(o,r)\backslash B(o,2) \}\Bigr)\\
    &\ge & \frac{a}{\gamma} \; \#\{L_i: L_i\cap L_{0,\varepsilon}\not=\emptyset\}
\end{eqnarray*}
and, since  $\#\{L_i: L_i\cap L_{0,\varepsilon}\not=\emptyset\}\sim
Poi(2\,\varepsilon\gamma/\pi)$, 
\[
    \liminf_{\gamma\rightarrow \infty}\E X_{\gamma,\xi}\;\ge\;    \lim_{\gamma\rightarrow \infty}  \frac{a}{\gamma} \; 
    \E\#\{L_i: L_i\cap L_{0,\varepsilon}\not=\emptyset\}
    \;=\; \frac{2\,\varepsilon a}{\pi} >0\,,
\]
which is a contradiction to (\ref{con.iks.zer}). Thus, $\xi=1$ holds. 
If the tessellation $T=T_I/T_{II}$ is a superposition/nesting such that $T_I$ is a PLT, then 
\begin{eqnarray*}
    \frac{1}{\gamma} \; \nu_1(\widetilde{T}_\gamma^{(1)}\backslash \widetilde{T}_{\gamma, \xi, \varepsilon}^{(1)}\cap B(o,r))
    &\ge & \ind_{\{o\in \widetilde{T}_{I,\gamma}^{(1)}\}}\; \frac{1}{\gamma}\; 
    \nu_1(\widetilde{T}_{I,\gamma}^{(1)}\backslash \widetilde{T}_{I,\gamma,\varepsilon}^{(1)}\cap B(o,r))\,,
\end{eqnarray*}
where $\widetilde{T}_{I,\gamma}^{(1)}$ denotes the part of $\widetilde{T}_\gamma^{(1)}$ which corresponds to $T_I$.
Since $T_I$ is assumed to be a PLT, the same arguments as  above can be applied to show that $\xi = 1$.
Finally, let $T$ be a PDT and let $N(y)$ denote that node of $T$ which is closest to $y\in \R^2$.
It has been
shown in \cite{Bac2000} that for any $t>0$ and $y\in\partial B(o,1)$, there is a path $P(ty)$  from  $N(o)$ to $N(ty)$ on $T^{(1)}$
with length $c(P(ty))$ such that almost
surely
\begin{equation}\label{bac.tou.res}
\lim_{t\longrightarrow \infty} \frac{c(P(ty))}{t}=\frac{4}{\pi}\; .
\end{equation}
 Consider the stationary point process $T^{(1)}\cap L$ of intersection points $\{X_i\}$, where $L=\{sy:\, s\in\R\}$ and
 $\cdots<X_{-1}<X_0\le 0 < X_1<\cdots$, 
 and denote by $c(X_i,X_j)$  the shortest path length from
$X_i$ to $X_j$ on $T^{(1)}$. Furthermore, consider the stationary marked point process $\{(X_i,c(X_i,N(X_i))\}$ and denote 
its typical mark by $c_N^*$.
For each $i>0$ we then have
\begin{eqnarray*}
    \frac{c(X_0,X_i)}{|X_i-X_0|}&\le& \frac{c(N(o),N(X_i))}{|X_i-X_0|}+\frac{c(X_0,N(o))}{|X_i-X_0|}+\frac{c(X_i,N(X_i))}{|X_i-X_0|}\\
    &\le & \frac{c(P(X_i))}{|X_i|}+\frac{c(X_0,N(o))}{|X_i|}+\frac{c(X_i,N(X_i))}{|X_i|}\;.
\end{eqnarray*}
Clearly,  the second summand of the latter expression tends to $0$ as $i\to\infty$. The same is true for the third summand,
because  $\{(X_i,c(X_i,N(X_i))\}$ is ergodic and
$\E \,c_N^*<\infty$.  
Thus, by (\ref{bac.tou.res}), we get that
\begin{equation}\label{lim.sup.lim}
\limsup\limits_{i\to\infty}  \frac{c(X_0,X_i)}{|X_i-X_0|}\le \frac{4}{\pi}\; .
\end{equation}
On the other hand, we have $\mathbb{P}(\lim_{i\to\infty}  c(X_0,X_i)/|X_i-X_0|=\xi)=1$ if and only if
$\mathbb{P}(\lim_{i\to\infty} c(X^*_0,X^*_i)/|X^*_i-X^*_0|=\xi)=1$, where $\{X_i^*\}$ is the Palm version of $\{X_i\}$.
Now, using (\ref{lim.sup.lim}) and (\ref{eq.xi}), it follows that $\xi\le 4/\pi$.
\end{proof}

\begin{figure}[tp]
\centering
\subfigure[$\kappa\to 0$ (where $\lambda_\ell=1$)]{
\includegraphics[width=7.0cm,viewport = 30 30 500 500,clip]{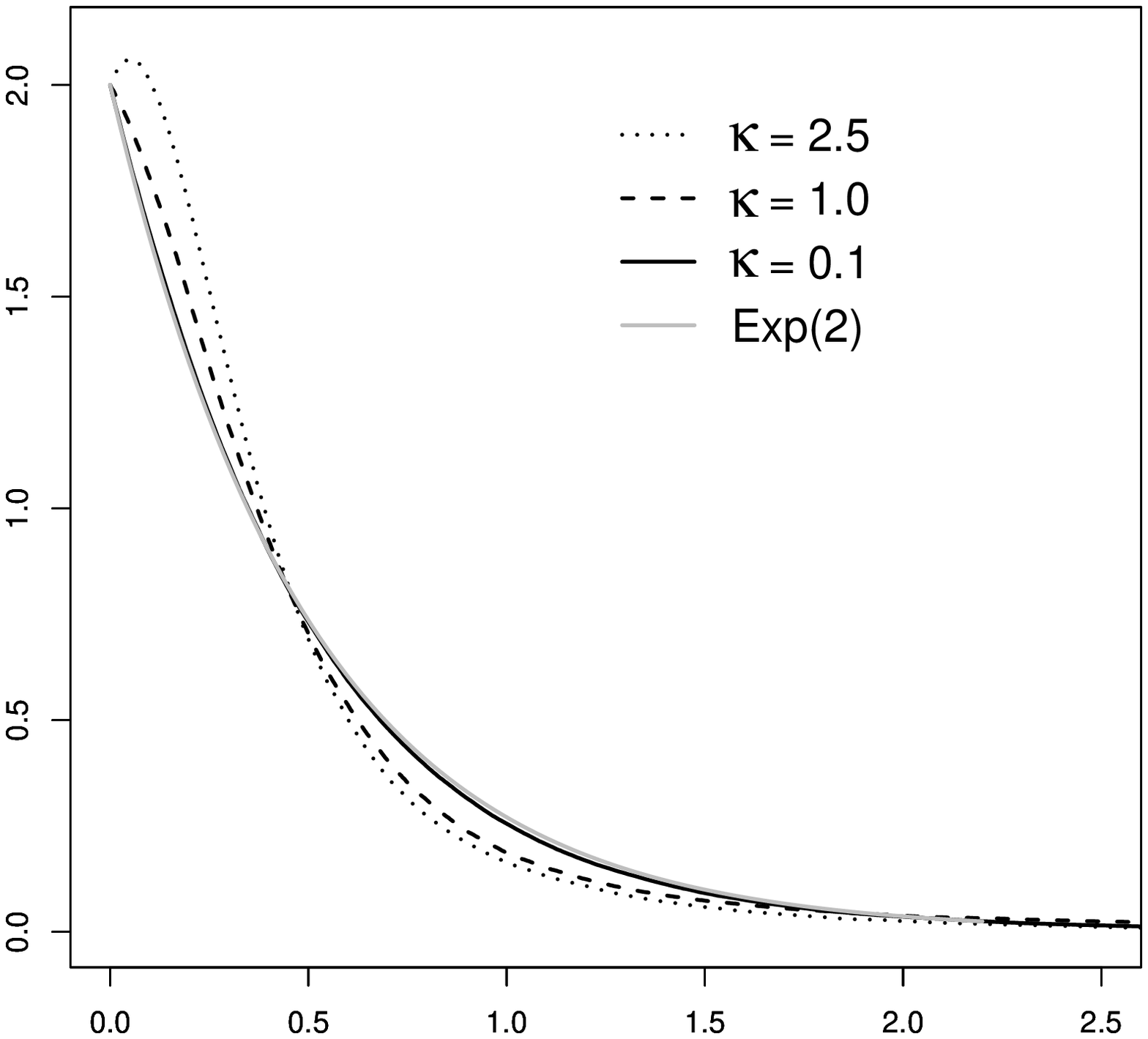}\label{fig:PLT.kappa.small}
}
\subfigure[$\kappa\to \infty$ (where $\gamma\lambda_\ell=1$)]{
\includegraphics[width=7.0cm,viewport = 30 30 500 500,clip]{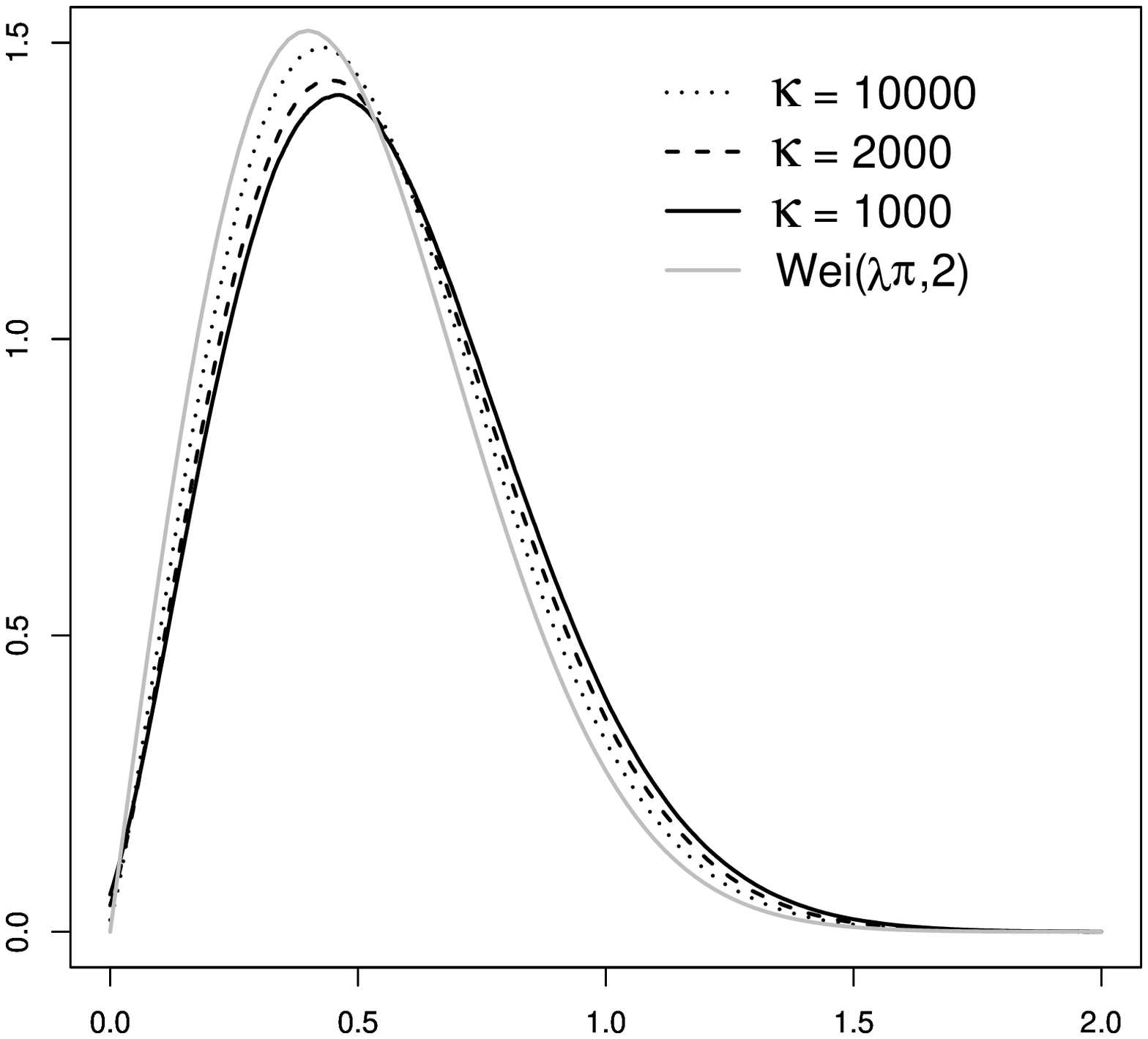}\label{fig:PLT.kappa.large}
}
\caption{Densities of $C^*$ if $T$ is a PLT (together with corresponding limit distributions)}\label{fig:PLT.densities}
\end{figure}

\subsection{Some extensions}\label{sub.fur.tes}

Note that the setting of Theorem~\ref{the.gam.inf} can be generalized in different ways.
For example, the statement of  this theorem remains valid if
instead of $C^*$  the typical subscriber line length $S^*$ is considered, where $S^*$ is the shortest path length
        from the origin to the nearest point $X_{H,0}$ of $X_H$, which is defined as the sum of the  distance from the origin
        to the nearest point of the edge set $T^{(1)}$ and the shortest path length on $T^{(1)}$ from this point to $X_{H,0}$
          (\cite{Glo07}). Note that in this case the auxiliary results corresponding to Lemmas~\ref{lem.exp.conv} 
          and \ref{lem.conv.prob} can be proved basically in the same way.

Furthermore, in the proof of Theorem~\ref{the.gam.inf} it is not necessary to assume that $T$ is a random tessellation, 
but it is possible to consider an arbitrary stationary and isotropic segment process in $\R^d$ which is mixing
        and such that there is only one single cluster
        with probability $1$. This means in particular that Theorem~\ref{the.gam.inf} can be extended to random geometric graphs.

Another kind of extensions can be obtained by relaxing the assumption         
     that $X_{L,n}$ is connected to the nearest point of $X_H$, i.e., $T_H$ is a Voronoi tessellation.
     For instance, $X_{L,n}$ can be connected to its $k$-th nearest neighbour of $X_H$ for any $k\ge 1$. 
         Then, in Theorem~\ref{the.gam.inf} we
         only have to replace $Z$ by the distance from the origin to the $k$-th nearest point of a Poisson process which
        is distributed according to a generalized Gamma distribution (\cite{Hae05,Zuy99}). Further possible extensions include
        that $T_H$ is a certain Cox-Laguerre tessellation (\cite{Lau08}) or an aggregated
        tessellation (\cite{Ba02,Tch01}).

\section{Conclusion and Outlook}\label{sec.con.out}

\begin{figure}
\subfigure[$\kappa\to 0$ (where $\lambda_\ell=1$)]{
\includegraphics[width=7.0cm,viewport = 30 30 500 500,clip]{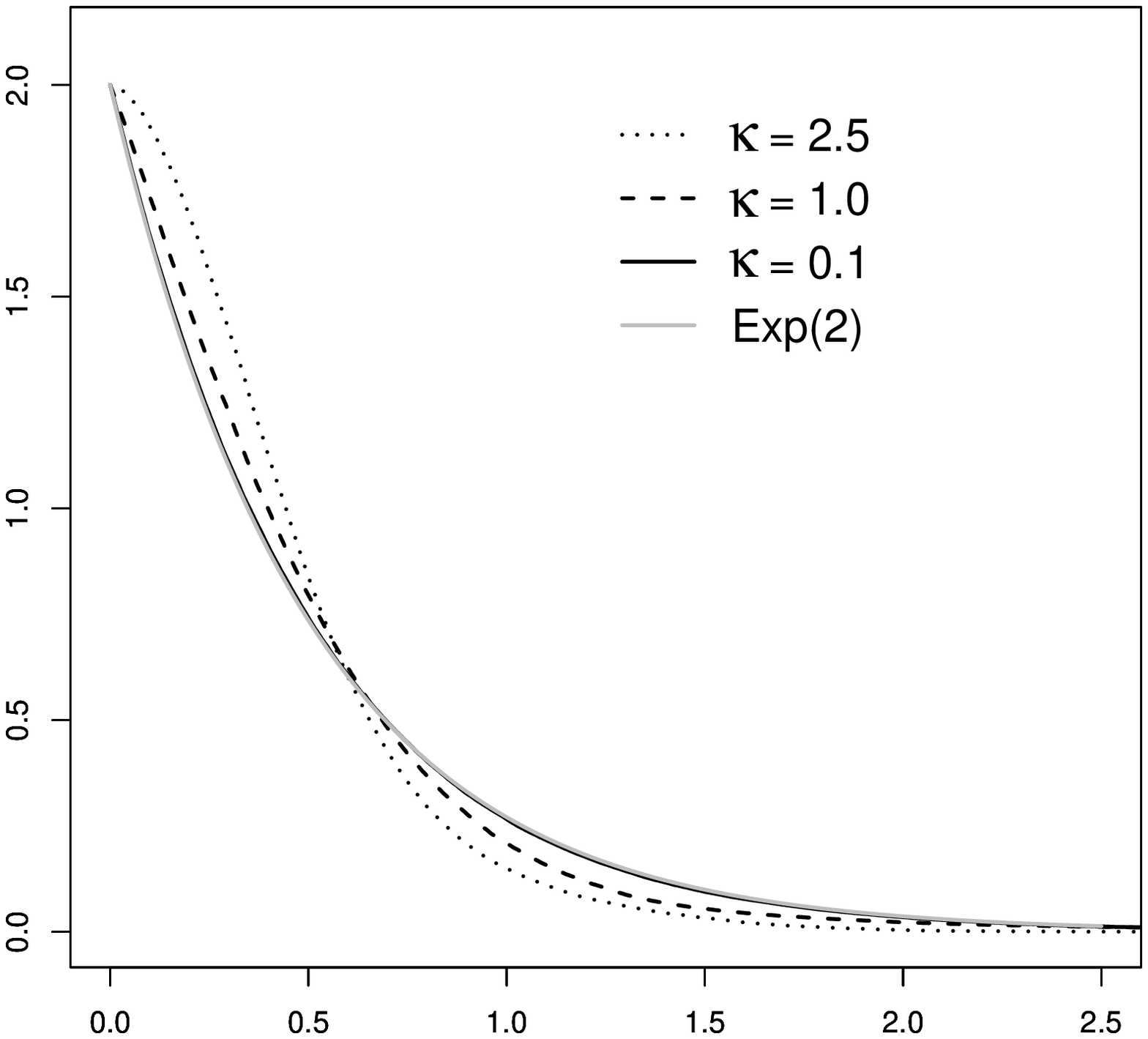}\label{fig:PVT.kappa.small}
}
\subfigure[$\kappa\to \infty$ (where $\gamma\lambda_\ell=1$)]{
\includegraphics[width=7.0cm,viewport = 30 30 500 500,clip]{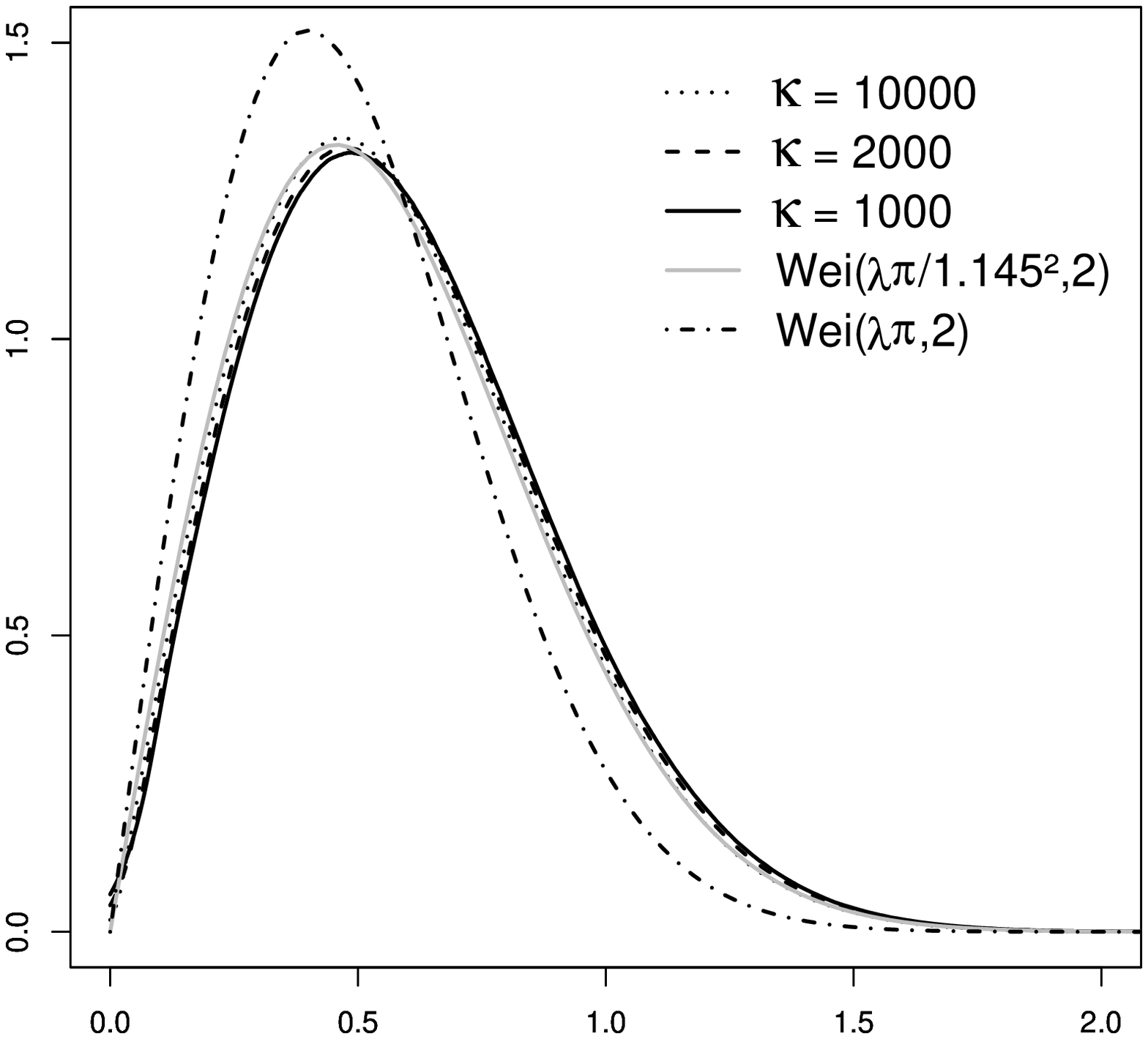}\label{fig:PVT.kappa.large}
}
\caption{Densities of $C^*$ if $T$ is a PVT (together with corresponding limit distributions)}\label{fig:PVT.densities}
\end{figure}
We consider the typical shortest path length $C^*$ of stochastic network models with two hierarchy levels, 
where the locations of network components are modelled by Cox processes on the edges of random tessellations. It is shown that the distribution
of $C^*$ converges to known limit distributions for extreme cases of the model parameters, i.e., if a certain scaling factor $\kappa$
tends to zero or infinity.  

The results of the present paper have applications in the analysis of telecommunication access networks since
the distribution of $C^*$ is closely related to cost and risk analysis of such networks (\cite{Glo09}).
Using the fitting techniques introduced in \cite{Glo06}, an optimal tessellation model can be chosen for a given set of road data.
Moreover, the scaling factor $\kappa$ can be estimated. Then, on the one hand,
for small values of $\kappa$ the limit distribution of $C^*$ is directly available and it does not depend on the type of
the optimal tessellation model.
On the other hand, for large values of $\kappa$ the limit distribution of $C^*$ and an upper bound for this distribution
is directly available if the optimal model is PLT or PLT-superposition/nesting and PDT, respectively.

In order to get an idea how small or large the scaling factor $\kappa$ should be (to replace the distribution of $C^*$
 by the corresponding limit distribution)  and how to calculate the constant $\xi$ appearing in the
limit distribution for $C^*$ as $\kappa\to\infty$, the density of $C^*$
can be estimated by Monte Carlo simulation of the typical serving zone (\cite{Voss08b}). This can be done for PVT, PLT and PDT
as well as for superpositions and nestings built from these basic tessellation models,
using simulation algorithms of the typical serving zone introduced in \cite{FL08,Glo05,Voss08a,Voss09b}.  
 In Figures~\ref{fig:PLT.densities} and~\ref{fig:PVT.densities}
estimated densities for different
values of $\kappa$ are shown together with the corresponding limit distributions if the tessellation model chosen for the underlying road system
is a PLT and PVT, respectively.
As can be seen in Figure~\ref{fig:PVT.densities}~(b),
the density of the $Wei(\lambda\pi/1.145^2,2)$-distribution approximates the density of $C^*$
very well for $T$ being a PVT and $\kappa \ge 1000$. This suggests that in this case the constant $\xi$
appearing in Theorem~\ref{the.gam.inf} and Corollary~\ref{cor.tes.mod}, respectively, is approximately $1.145$.
\begin{figure}
\centering
\subfigure[PLT]{
\includegraphics[width=7.0cm,viewport = 30 30 500 500,clip]{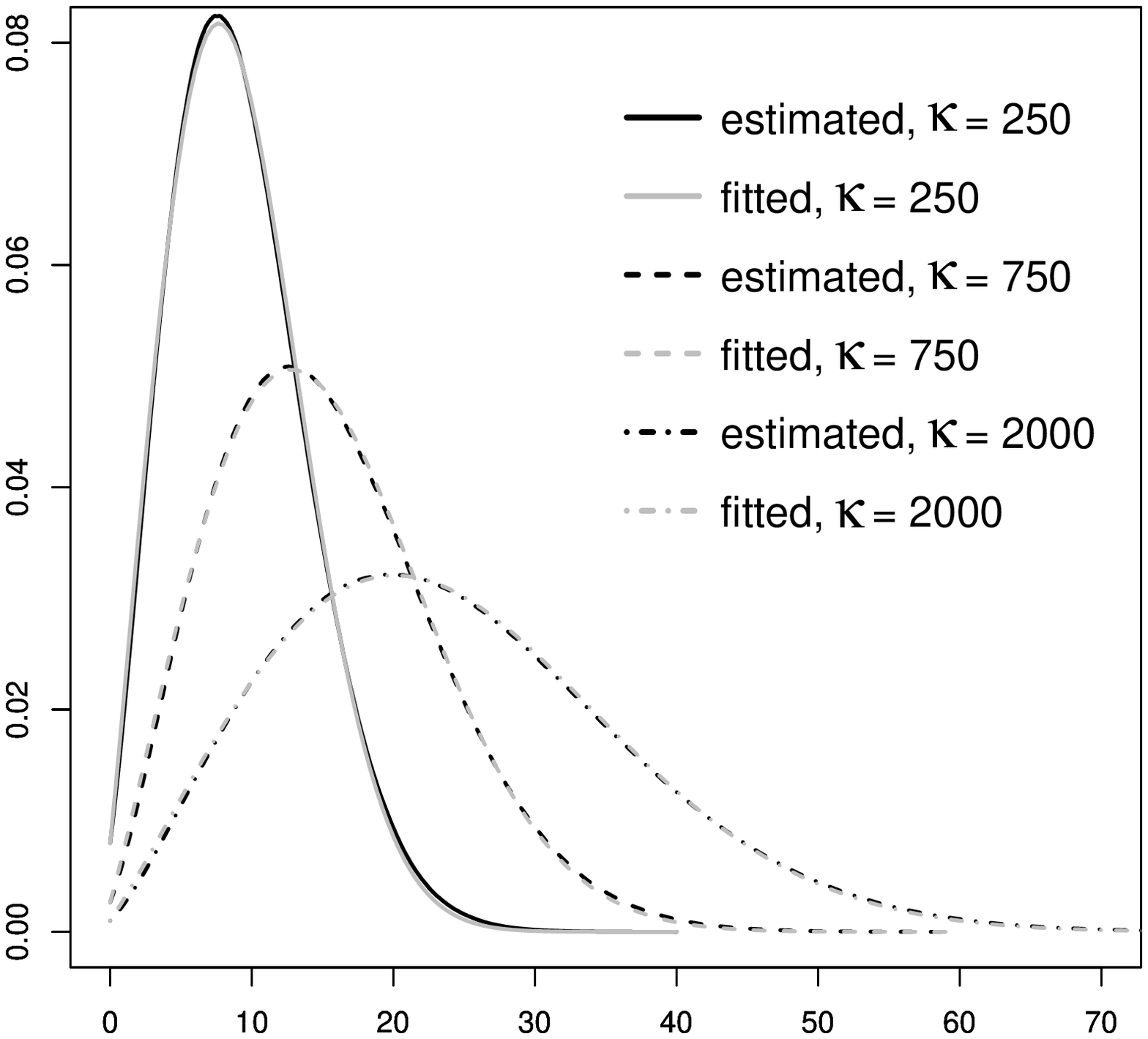}
}
\subfigure[PVT]{
\includegraphics[width=7.0cm,viewport = 30 30 500 500,clip]{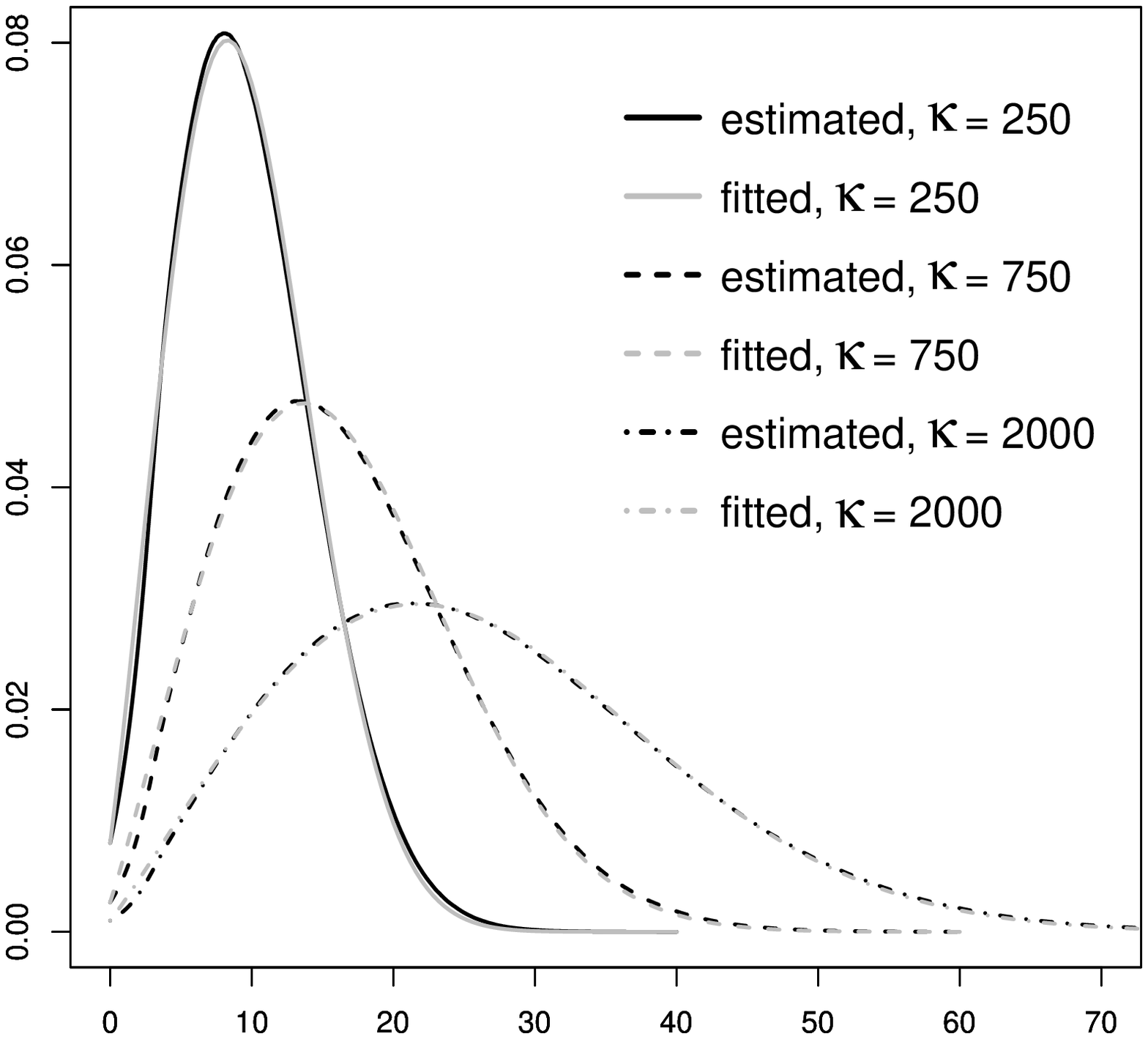}
}
\caption{Densities of $C^*$ (together with fitted parametric densities)}\label{fig:fitted.density}
\end{figure}

Furthermore, the limiting distributions derived in the present paper can be used to choose parametric
densities which can be fitted to the estimated density of $C^*$ for a large range of $\kappa$. Parametric families
which include both exponential distributions and Weibull distributions turned out to be
good choices, see \cite{Glo09}. In Figure~\ref{fig:fitted.density} estimated densities for different values of $\kappa$
are shown together with  fitted truncated Weibull distributions. Note that these
truncated Weibull distributions have two parameters and there is a quite good fit for both tessellation 
models considered in Figure~\ref{fig:fitted.density} and for a large range of values of $\kappa$.




\begin{appendix}

\vspace{0.5cm}

\section{Some mathematical background}\label{app.ein.ein}

In the proof of Lemma~\ref{lem.exp.conv} given below we make use
of some well-known results from measure theory, the theory of
subadditive processes, and geometric measure theory which are
briefly summarized. We start with the definition of convergence in
measure and uniform integrability which can be used to
characterize $L^1$-convergence. A family of measurable functions
$\{f_\gamma,\gamma\ge 1\}$ defined on a measurable space
$(\Omega,\mathcal{A},\mu)$ and
taking values in $\R$ converges locally
in $\mu$-measure to a measurable function $f:\Omega\to\R$ if
\begin{equation}\label{def.con.mea}
    \lim_{\gamma\rightarrow \infty}\mu(\{|f_\gamma-f|\ge \varepsilon\}\cap A)=0
\end{equation}
for all $\varepsilon>0$ and  $A\in\mathcal{A}$ with
$\mu(A)<\infty$,  where $\mu$ is assumed to be a $\sigma$-finite measure.   If $\mu$ is a probability measure such that
(\ref{def.con.mea}) holds for each $\varepsilon>0$ and $A=\Omega$,
then one says that $f_\gamma$ converges in probability to $f$.
Furthermore, if for each $\varepsilon > 0$ there is a
$\mu$-integrable function $g$ such that
\begin{equation}\label{eq.uni.int}
    \int_{\{|f_\gamma| \ge g\}} |f_\gamma(\omega)|\mu(d\omega) \le \varepsilon\quad\mbox{for all }\gamma\ge 1\,,
\end{equation}
then the family $\{f_\gamma,\gamma\ge 1\}$ is said to be uniformly
$\mu$-integrable. With the above definitions it is possible to
characterize the $L^1$-convergence as follows; see Theorem~2.12.4
in \cite{Bauer81}.
\begin{thm}\label{lem.uni.int}
    A sequence of $\mu$-integrable functions $f_1,f_2,\dots:\Omega\to\R$ converges in $L^1$ to a $\mu$-integrable function $f:\Omega\to\R$
    if and only if {\rm (i)} $f_n$ converges locally in $\mu$-measure to $f$ and {\rm (ii)} $\{f_n\}$ is
uniformly $\mu$-integrable.
\end{thm}
We still mention an elementary but useful result which immediately
follows from the definition of uniform integrability.
\begin{lem}\label{lem.dom.int}
    Let $\{f_\gamma,\gamma\ge 1\}$ and $\{g_\gamma,\gamma\ge 1\}$ be two families of measurable functions on $(\Omega,\mathcal{A},\mu)$
    which satisfy that
    $|f_\gamma| \le |g_\gamma|$ for all $\gamma\ge 1$. Then
    $\{f_\gamma,\gamma\ge 1\}$ is uniformly $\mu$-integrable if $\{g_\gamma,\gamma\ge 1\}$ is uniformly $\mu$-integrable.
\end{lem}
Another useful tool is the notion of subadditivity. Let
$\mathbf{Y}$ be a family of real-valued random variables
$\mathbf{Y}=\{Y_{ij}, 0\le i<j\}$ defined on some probability
space $(\Omega,\mathcal{A},\mathbb{P})$. Note that $\mathbf{Y}$
can be seen as a random element of some measurable space
$(\mathcal{S},\mathcal{B}(\mathcal{S}))$ of double-indexed
sequences, where $\mathcal{B}(\mathcal{S})$ is the
Borel-$\sigma$-algebra of $\mathcal{S}$. Then $\mathbf{Y}$
 is called a subadditive process
    if
    \begin{enumerate}
        \item $Y_{ik}\le Y_{ij}+Y_{jk}$ for all $i < j < k$,
        \item $\mathbf{Y}=\{Y_{ij}\}\stackrel{{\rm d}}{=}\mathbf{Y}'=\{Y_{i+1,j+1}\}$,
        \item $\E Y_{01}^+ <\infty$, where
        $Y_{01}^+=\max\{0,Y_{01}\}$.
    \end{enumerate}
The following result is due to Kingman (\cite{Kingman73},
Theorem~1). It is called the subadditive ergodic theorem; see also
Theorem~10.22 in \cite{Kal02}.
\begin{thm}\label{th.sub.erg}
    Let $\mathbf{Y}$ be a subadditive process. Then the limit
    \begin{equation}\label{con.sub.add}
        \zeta=\lim_{j\rightarrow \infty}\frac{1}{j}Y_{0j}
    \end{equation}
    exists and is finite with probability one and $\E \zeta=\inf_{j\in\N}\E Y_{0j}/j$.
    If $\E \zeta > -\infty$, then the convergence in $(\ref{con.sub.add})$ also holds in the  $L^1$-norm.
    Moreover, let $\mathcal{I}_\mathcal{S}\subset\mathcal{B}(\mathcal{S})$ be the  $\sigma$-algebra of subsets of $\mathcal{S}$
    which are invariant under the shift
    $\mathbf{Y}\mapsto \mathbf{Y}'$, where $Y_{ij}'=Y_{i+1,j+1}$,
    and let $\mathcal{I}=\mathbf{Y}^{-1}\mathcal{I}_\mathcal{S}\subset\mathcal{A}$
    be the corresponding sub-$\sigma$-algebra of events.
    Then,
    \begin{equation}\label{zet.erg.cas}
        \zeta=\lim_{j\rightarrow \infty}\frac{1}{j}\;\E\big(Y_{0j}\mid\mathcal{I}\big)\,.
    \end{equation}
\end{thm}

Note that a subadditive process $\mathbf{Y}$ is called ergodic if
$\mathbb{P}(\mathbf{Y}\in A)=0$ or $\mathbb{P}(\mathbf{Y}\in A)=1$
for each $A\in\mathcal{I}_\mathcal{S}$. Thus, in the ergodic case,
the limit $\zeta$ considered in (\ref{con.sub.add}) and
(\ref{zet.erg.cas}), respectively, is almost surely constant.

Finally, we use a decomposition of the Hausdorff measure $\nu_1$
which is a special case of the generalized Blaschke-Petkantschin
formula (\cite{jen98}, Proposition 5.4).
\begin{thm}\label{thm.blaschke}
    Let $C\subset\R^2$ be a differentiable curve and assume that
    \begin{equation}\label{nuu.set.nul}
        \nu_1(\{x\in C:\,\mbox{Tan}[C,x]= \mbox{span}\{x\}\})=0\,,
    \end{equation}
    where $\mbox{Tan}[C,x]$ is the tangent at $x$ to $C$ and $\mbox{span}\{x\}=\{cx:c\in\R\}$ is the line which goes through the origin
    $o\in\R^2$ and  the point $x\in C$.
    Then, for any measurable $g:C\to [0,\infty)$ it holds that
    \begin{equation}\label{lem.blaschke.eq}
        \int_{C}g(x)\,\nu_1(dx)=\int_0^{2\pi} \sum_{x_i\in C\cap  L^+_\Phi} \frac{|x_i|}{\sin\alpha_i}\;g(x_i)\;d\Phi\,,
    \end{equation}
    where $L^+_\Phi$ is the half line
    of direction $\Phi\in[0,2\pi)$ emanating from $o$ and $\alpha_i$ is the angle between $\mbox{Tan}[C,x_i]$ and $\mbox{span}\{x_i\}$.
\end{thm}

\section{Proof of Lemma~\ref{lem.exp.conv}}\label{app.zwe.zwe}

With the help of Theorems~\ref{lem.uni.int} -- \ref{thm.blaschke}
stated above we are now able to prove Lemma~\ref{lem.exp.conv}.
Obviously, $        \limsup_{\gamma\rightarrow \infty}\E
        \exp\Big(-\frac{\lambda}{\gamma}\;
\nu_1\big(\widetilde{T}^{(1)}_\gamma\setminus\widetilde{T}^{(1)}_{\gamma, \xi, \varepsilon}\cap
B(o,r)\big)\Big)\le 1$.
    Thus it is sufficient to show that
    \begin{equation}\label{lim.inf.gam}
        \liminf_{\gamma\rightarrow \infty}\E
        \exp\Big(-\frac{\lambda}{\gamma}\;
\nu_1\big(\widetilde{T}^{(1)}_\gamma\setminus\widetilde{T}^{(1)}_{\gamma, \xi, \varepsilon}\cap
B(o,r)\big)\Big)\ge 1\,.
    \end{equation}

\vspace{0.5cm} \noindent
 {\bf Proof of (\ref{lim.inf.gam}).} $\;$   First recall that we can identify $\widetilde{T}^{(1)}_\gamma$ with the Palm version
    $\Lambda^*_{T^{(1)}_\gamma}$ of the
    stationary random measure $\Lambda_{T^{(1)}_\gamma}$ given by
    $\Lambda_{T^{(1)}_\gamma}(B)=\nu_1(B\cap T^{(1)}_\gamma)$ for $B\in \mathcal{B}^2$ since
    $\Lambda_{T^{(1)}_\gamma}$ is the random driving measure of the Cox process $X_L$, see \cite{Stoyan95}, p. 156.
    Then, using the abbreviation
\[
    h(\tau^{(1)})=\exp\Bigl(-\frac{\lambda}{\gamma}\,
\nu_1(\tau^{(1)}\setminus \tau^{(1)}_{\xi,\varepsilon}\cap
B(o,r))\Bigr)\,,
\]
 where
$\tau^{(1)}_{\xi,\varepsilon}=\big\{u\in \tau^{(1)}:\,
\big|c(u)-\xi|u|\big| < \varepsilon\big\}$ and $c(u)$ denotes the
length of the shortest path from $u$ to the origin along the edge
set $\tau^{(1)}$ of a tessellation $\tau$ with $o\in \tau^{(1)}$,
we
    get from the Campbell theorem for stationary random measures (\cite{Daley0307}, Proposition~13.2.V) that
    \begin{eqnarray*}
        \E h(\widetilde{T}^{(1)}_\gamma)&=&\frac{1}{\gamma\nu_2(B(o,1/\gamma))}
        \E\Big(\int_{T^{(1)}_\gamma\cap B(o,1/\gamma)}h(T^{(1)}_\gamma -x )\nu_1(dx)\Big)\\
        &=&\frac{1}{\pi}
        \E\Big(\int_{T^{(1)}\cap B(o,1)}h\big(T^{(1)}_\gamma -\frac{z}{\gamma} \big)\nu_1(dz)\Big)\,,
    \end{eqnarray*}
    where we used the substitution $z=\gamma x$ in the last expression bearing in mind that
    $(1/\gamma) T^{(1)}=T^{(1)}_\gamma$. Furthermore, we put $T^{(1)}_{\gamma,\varepsilon,z}=\{y\in T^{(1)}_\gamma:
    |c(y,z/\gamma)-\xi|y-z/\gamma||< \varepsilon\}$, where
    $c(y,z/\gamma)$ denotes the length of the shortest path from $y$ to
    $z/\gamma$ along the edges of the considered graph. Then, for each $\gamma\ge 1$,  we
    get that
    \begin{eqnarray*}
        \lefteqn{\E h\big(\widetilde{T}_\gamma^{(1)}\big)
        =\frac{1}{\pi}\;
        \E\Big(\int_{T^{(1)}\cap B(o,1)}\exp
        \Big(-\frac{\lambda}{\gamma}\;
        \nu_1 \big(T^{(1)}_\gamma\backslash T^{(1)}_{\gamma,\varepsilon,z}\cap B(z/\gamma,r )\big)\Big)\nu_1(dz)\Big)}\\
        &&\ge\frac{1}{\pi}\;
        \E\Big(\nu_1\big(T^{(1)}\cap B(o,1)\big)\inf_{z\in T^{(1)}\cap B(o,1)}\exp
        \Big(-\frac{\lambda}{\gamma}\;
        \nu_1 \big(T^{(1)}_\gamma\backslash T^{(1)}_{\gamma,\varepsilon,z}\cap B(z/\gamma,r )\big)\Big)\Big)\\
        &&=\frac{1}{\pi}\;
        \E\Big(\nu_1\big(T^{(1)}\cap B(o,1)\big)\exp \Big(-\sup_{z\in T^{(1)}\cap
        B(o,1)}\frac{\lambda}{\gamma}\;
        \nu_1 \big(T^{(1)}_\gamma\backslash T^{(1)}_{\gamma,\varepsilon,z}\cap B(z/\gamma,r )\big)\Big)\Big)\\
        &&\ge\frac{1}{\pi}\;
        \E\Big(\nu_1\big(T^{(1)}\cap B(o,1)\big)\exp \Big(-\sup_{z\in T^{(1)}\cap
        B(o,1)}\frac{\lambda}{\gamma}\;
        \nu_1 \big(T^{(1)}_\gamma\backslash T^{(1)}_{\gamma,\varepsilon,z}\cap B(o,r + 1)\big)\Big)\Big)\,.
    \end{eqnarray*}
    Now, in order to prove (\ref{lim.inf.gam}),
    it is sufficient to show that
 \begin{equation}\label{eq.to.show1}
    X_{\gamma,\xi} \stackrel{L^1}{\rightarrow} 0\qquad\mbox{ for
    $\gamma\to\infty$,}
 \end{equation}
    where $
        X_{\gamma,\xi}=\sup_{z\in T^{(1)}\cap
        B(o,1)}\frac{1}{\gamma}\;
        \nu_1 \Big(T^{(1)}_\gamma\backslash T^{(1)}_{\gamma,\varepsilon,z}\cap B(o,r +
        1)\Big)$.
    To see this, note first that (\ref{eq.to.show1}) implies that $X_{\gamma,\xi}$ converges in probability to $0$.
    Thus,  the random variable $Y_{\gamma,\xi}=\exp(-\lambda
X_{\gamma,\xi})\nu_1(T^{(1)}\cap B(o,1))$ converges in probability to
    $\nu_1\big(T^{(1)}\cap B(o,1)\big)$ if (\ref{eq.to.show1}) holds. Moreover,
    $Y_{\gamma,\xi}\le \nu_1(T^{(1)}\cap B(o,1))$ for all $\gamma \ge 1$ and $\E\nu_1(T^{(1)}\cap B(o,1))=\pi
    <\infty$, which means that $\{Y_{\gamma,\xi},\,\gamma\ge 1\}$ is uniformly integrable. Hence,
    Theorem~\ref{lem.uni.int} yields that $Y_{\gamma,\xi}$ converges in $L^1$ to $\nu_1\big(T^{(1)}\cap
    B(o,1)\big)$ and, in particular, $\lim_{\gamma\rightarrow \infty}1/\pi\E Y_{\gamma,\xi}=1/\pi \E \nu_1\big(T^{(1)}\cap B(o,1)\big)=1$
    if (\ref{eq.to.show1}) holds.
    Thus, (\ref{lim.inf.gam}) follows if we can show that (\ref{eq.to.show1})  is true.

\vspace{0.5cm} \noindent
 {\bf Proof of (\ref{eq.to.show1}).} $\;$
Since $X_{\gamma,\xi}\ge 0$ it suffices to show that $\E _{\gamma,\xi}\to 0$.
Furthermore, note that with probability $1$ the segments of the
segment system $T_\gamma^{(1)}\cap B(o,r+1)$ fulfill the conditions
of Theorem~\ref{thm.blaschke}, since none of these segments
,,points'' to the origin because $T^{(1)}$ was assumed to be stationary. Thus, using Theorem~\ref{thm.blaschke}
we get that
    \begin{eqnarray*}
        \E X_{\gamma,\xi}\!\!\!\! &=&\!\!\! \E\Big(\sup_{z\in T^{(1)}\cap B(o,1)}\frac{1}{\gamma}
        \int_{T_\gamma^{(1)}\cap B(o,r+1)}\ind_{[\varepsilon,\infty)}\Big(\big|c(y,\frac{z}{\gamma})-\xi|y-\frac{z}{\gamma}|\big|\Big)\,\nu_1(dy)\Big)\\
        &=&\!\!\!\E\Big(\frac{1}{\gamma}\!\sup_{z\in T^{(1)}\cap B(o,1)}
        \int_0^{2\pi}\!\!\!\! \sum_{\substack{X_i\in T^{(1)}_\gamma\cap L^+_\Phi:\\|X_i|\le r+1}}\,\frac{|X_i|}{\sin\alpha_i}
        \ind_{[\varepsilon,\infty)}\Big(\big|c(X_i,\frac{z}{\gamma})-\xi|X_i-\frac{z}{\gamma}|\big|\Big)\,d\Phi\Big)\\
        &\le&\!\!\!\!\! \frac{r+1}{\gamma}\E\Big(\!\!\int_0^{2\pi}\!\!\!\!
        \sup_{z\in T^{(1)}\cap B(o,1)} \sum_{
        \substack{X_i\in T^{(1)}_\gamma\cap L^+_\Phi:\\|X_i|\le r+1}}\!\!\!\frac{1}{\sin\alpha_i}
        \ind_{[\varepsilon,\infty)}\Big(\big|c(X_i,\frac{z}{\gamma})-\xi|X_i-\frac{z}{\gamma}|\big|\Big)\,d\Phi\!\!\Big)\\
        &=&\!\!\!\frac{2\pi(r+1)}{\gamma}\E\Big(\!
        \sup_{z\in T^{(1)}\cap B(o,1)} \sum_{
        \substack{X_i\in T^{(1)}_\gamma\cap L^+:\\|X_i|\le r+1}}\!\!\!\frac{1}{\sin\alpha_i}
        \ind_{[\varepsilon,\infty)}\Big(\big|c(X_i,\frac{z}{\gamma})-\xi|X_i-\frac{z}{\gamma}|\big|\Big)\Big)\\
&=&\!\!\! 2\pi(r+1)   \E g_\gamma\big(T^{(1)}\big)
        \,,
    \end{eqnarray*}
    where in the last but one line we used Fubini's theorem and the isotropy of $T^{(1)}_\gamma$, denoting by $L^+=L^+_0$
    the half line with direction $\Phi=0$, and in the last expression we
     used the abbreviation
    \begin{equation}\label{def.fun.gam}
        g_\gamma\big(T^{(1)}\big)=\frac{1}{\gamma}\sup_{z\in T^{(1)}\cap B(o,1)} 
        \sum_{\substack{X_i\in T^{(1)}_\gamma\cap L^+:\\|X_i|\le r+1}}\,\frac{1}{\sin\alpha_i}
        \ind_{[\varepsilon,\infty)}\Big(\big|c(X_i,\frac{z}{\gamma})-\xi|X_i-\frac{z}{\gamma}|\big|\Big)\,.
    \end{equation}
    Since the point process $T^{(1)}\cap \R$  is stationary with intensity $2/\pi$ (\cite{ScW08}, Theorem 4.5.3), where we identify
    $\R$ with the x-axis,
    we can apply the inversion formula for Palm distributions of stationary point processes on $\R$; see Proposition 11.3 ~(iii) in \cite{Kal02}.
    Thus, if $T^{(1)*}$ denotes the Palm version of $T^{(1)}$ with respect to the point process $T^{(1)}\cap \R$,
    we get that
    \begin{eqnarray*}
        \E g_\gamma\big(T^{(1)}\big) &=& \frac{2}{\pi}\;\E\Big(\int_0^\infty  \ind_{[0,X_1^*]}(x)\,g_\gamma\big(T^{(1)*}-x\big)\, dx\Big)\,,
    \end{eqnarray*}
    where the points of $\{X_i^*\}=T^{(1)*}\cap \R$ are numbered
    in ascending order
    such that  $\ldots<X_{-1}^*<X_0^*= 0 < X_1^*<X_2^*<\ldots$. Hence, in order to prove (\ref{eq.to.show1}) it suffices to show that
    \begin{equation}\label{eq.to.show2}
\lim\limits_{\gamma\to\infty}\E\Big(\int_0^\infty
\ind_{[0,X_1^*]}(x)\,g_\gamma\big(T^{(1)*}-x\big)\, dx\Big)=0\,,
\end{equation}
where the function $g_\gamma:{\cal F}\to [0,\infty)$ is given in (\ref{def.fun.gam}).
The proof of (\ref{eq.to.show2}) is subdivided into two main steps. First, we
 show that
 \begin{equation}\label{alm.shu.con}
 \lim_{\gamma\to\infty} \widetilde{g}_\gamma\big(x,T^{(1)*}\big)=0
 \end{equation}
 almost everywhere with respect to the product measure $\nu_1\otimes \mathbb{P}^*$, where we used the abbreviating
 notation $\widetilde{g}_\gamma(x,T^{(1)*})=\ind_{[0,X_1^*]}(x)\,g_\gamma(T^{(1)*}-x)$
 and $\mathbb{P}^*$ denotes the distribution of $T^{(1)*}$.
 Then, we show that $\{\widetilde{g}_\gamma,\,\gamma>0\}$ is uniformly $(\nu_1\otimes \mathbb{P}^*)$-integrable.
 By means of Theorem~\ref{lem.uni.int}, this implies that (\ref{eq.to.show2}) holds.

 \vspace{0.5cm} \noindent
 {\bf Proof of (\ref{alm.shu.con}).} $\;$
    Note that for each $x\in [0, X_1^*] $ we get
    \begin{eqnarray*}
        \lefteqn{ g_\gamma\big(T^{(1)*}-x\big)}\\
        &&\le\frac{1}{\gamma}\sup_{z\in (T^{(1)*}-x)\cap B(o,1)} 
        \sum_{\substack{X_i\in
        (T^{(1)*}_\gamma\!-\frac{x}{\gamma})\cap L^+:\\|X_i|\le r+1}}\!\!\!\frac{1}{\sin\alpha_i}
        \ind_{[\varepsilon,\infty)}\!\Big(\!\big|c(X_i,\frac{z}{\gamma})-\xi|X_i-\frac{z}{\gamma}|\big|\Big)\\
        &&=\frac{1}{\gamma}\sup_{z\in T^{(1)*}\cap B(x,1)} 
        \sum_{\substack{X_i^*\in
        T^{(1)*}\cap (L^++x):\\X_i^*\in B(x,(r+1)\gamma)}}\!\frac{1}{\sin\alpha_i}
        \ind_{[\varepsilon,\infty)}\Big(\frac{1}{\gamma}\big|c(X_i^*,z)-\xi|X_i^*-z|\big|\Big)\\
        &&\le \frac{1}{\gamma}\sum_{\substack{X_i^*\in
        T^{(1)*}\cap (L^++x):\\X_i^*\in B(x,(r + 1)\gamma)}}\!\!\frac{1}{\sin\alpha_i}
        \sup_{z\in T^{(1)*}\cap B(x,1)}\ind_{[\varepsilon,\infty)}\Big(\frac{1}{\gamma}\big|c(X_i^*,z)-\xi|X_i^*-z|\big|\Big)\,.
          \end{eqnarray*}
Thus,
         \begin{eqnarray*}
        \lefteqn{ g_\gamma\big(T^{(1)*}-x\big)}\\
        &&\le \frac{1}{\gamma}\sum_{\substack{X_i^*\in
        T^{(1)*}\cap L^+:\\|X_i^*|\le (r + a)\gamma}}\!\!\frac{1}{\sin\alpha_i}
        \sup_{z\in T^{(1)*}\cap B(o,a)}\ind_{[\varepsilon,\infty)}\Big(\frac{1}{\gamma}\big|c(X_i^*,z)-\xi|X_i^*-z|\big|\Big)\\
        &&= \frac{1}{\gamma}\sum_{\substack{X_i^*\in
        T^{(1)*}\cap L^+:\\|X_i^*|\le (r + a)\gamma}}\!\!\frac{1}{\sin\alpha_i}
        \ind_{[\varepsilon,\infty)}\Big(\frac{1}{\gamma}\;\sup_{z\in T^{(1)*}\cap B(o,a)}\big|c(X_i^*,z)-\xi|X_i^*-z|\big|\Big)
        \,,
    \end{eqnarray*}
    where $a=1+X_1^*$.  Furthermore, we have
    \begin{eqnarray*}
       \frac{1}{\gamma}\;\sup_{z\in T^{(1)*}\cap B(o,a)}\big|c(X_i^*,z)-\xi|X_i^*-z|\big|
&\le&\frac{1}{\gamma}\;\big(c(o,X_i^*)-\xi |X_i^*|\big)\nonumber\\
        &+&\frac{1}{\gamma}\;\Bigg(\sup_{z\in T^{(1)*}\cap B(o,a)}c(z,o)+\xi a\Bigg)\,,
    \end{eqnarray*}
    since $c(X_i^*,o)-c(o,z)\!\le\! c(X_i^*,z)\!\le\! c(X_i^*,o)+c(o,z)$ and $\xi |X_i^*|-\xi a\le \xi|X_i^*-z|\! \le\! \xi |X_i^*|+\xi a$
for all $i\ge 1$ and $z\in T^{(1)*}\cap B(o,a)$. Clearly, the
second term of this upper bound tends to zero
$\mathbb{P}^*$-almost surely as $\gamma\to\infty$.  Thus in order
to show that (\ref{alm.shu.con}) holds, it suffices to prove that
$\mathbb{P}^*$-almost surely
   \begin{equation}\label{equ.con.zer}
        \frac{1}{\gamma}\;\bigl( c(o,X_i^*)-\xi X_i^*\bigr)\in \Bigl(-\;\frac{\varepsilon}{2}\;,\;\frac{\varepsilon}{2}\Bigr)
    \end{equation}
    for all sufficiently large $i\ge 1$ such that $X_i^*\le (r + a)\gamma$.

    \vspace{0.5cm} \noindent
 {\bf Proof of (\ref{equ.con.zer}).} $\;$
Note that $\mathbf{X}=\{|X_i^*-X_j^*|,i,j\ge 1,i<j\}$ is an
additive process, because $|X_i^*-X_k^*|=
|X_i^*-X_j^*|+|X_j^*-X_k^*|$ for $i<j<k$. Since $T^{(1)*}\cap
\R$ is cycle-stationary (see e.g. \cite{Thor2000}), we have that
$\{|X_i^*-X_j^*|\}\stackrel{{\rm d}}{=}\{|X_{i+1}^*-X_{j+1}^*|\}$,
where $0< \E X_1^* <\infty$. Thus, by Theorem \ref{th.sub.erg} we
get that the finite limit $\lim_{i\to\infty}
X_i^*/i=\zeta_{\mathbf{X}}$ exists $\mathbb{P}^*$-almost surely.
Furthermore,
 consider the family
 $\mathbf{Y}=\{Y_{ij},i,j\ge 1,i<j\}$ of non-negative random variables with  $Y_{ij}=c(X_i^*,X_j^*)$, where $c(X_i^*,X_j^*)$
denotes the shortest path length
    from $X_i^*$ to $X_j^*$ on $T^{(1)*}$. Then,
    it is easy to see that $Y_{ik}\le Y_{ij}+Y_{jk}$ for $i<j<k$.
    By the cycle-stationarity of $T^{(1)*}\cap \R$, we have that
    $\{Y_{ij}\}\stackrel{{\rm d}}{=}\{Y_{i+1,j+1}\}$, where
$\E Y_{01} =\E c(X_0^*,X_1^*)<\infty$ holds by condition (\ref{int.con.deb}); see the next paragraph below. Thus
$\mathbf{Y}$ is a subadditive process and we
    can again apply  Theorem \ref{th.sub.erg} to get that the finite limit $\lim_{j\to\infty}c(X_0^*,X_j^*)/j=\zeta_{\mathbf{Y}}$ exists
$\mathbb{P}^*$-almost surely.
    Since $\mathbf{X}$ and $\mathbf{Y}$ are ergodic (see the paragraphs below), the limits $\zeta_{\mathbf{X}}$ and $\zeta_{\mathbf{Y}}$
are  constant. Noticing that $0<\E X_1^*=\zeta_{\mathbf{X}}\le
\zeta_{\mathbf{Y}}<\infty$, this gives that
\begin{equation}\label{eq.xi}
        \lim_{j\to \infty}\frac{c(o,X_j^*)}{X_j^*}=\lim_{j\to
        \infty}\frac{j}{X_j^*}\;\frac{c(X_0^*,X_j^*)}{j}=\xi\,,
    \end{equation}
where $\xi=\zeta_{\mathbf{Y}}/\zeta_{\mathbf{X}}\in[1,\infty)$.
Now let $\widetilde{\varepsilon}>0$
    such that $\widetilde{\varepsilon} (r+a) < \varepsilon/2$. Then (\ref{eq.xi}) implies
    that with probability $1$
    \begin{eqnarray*}
        \frac{c(o,X_i^*)}{X_i^*}-\xi\in (-\widetilde{\varepsilon},\widetilde{\varepsilon})
    \end{eqnarray*}
    for all $i$ sufficiently large and, therefore,
    \[
\frac{1}{\gamma}\;\bigl( c(o,X_i^*)-\xi X_i^*\bigr)
    \in \Bigl(-\;\frac{\varepsilon}{2}\;,\;\frac{\varepsilon}{2}\Bigr)
    \]
    if $i$ is sufficiently large  and $X_i^*/\gamma\le r + a$.

   \vspace{0.5cm} \noindent
 {\bf Proof of $\;\E c(X_0^*, X_1^*)\!<\!\infty$.} 
Consider the stationary marked point process $\{(X_n,\Xi_n^+)\}$,
where $\{X_n\}=T^{(1)}\cap \R$ is the point process of
intersection points of the edge set $T^{(1)}$ with the line $\R$,
and  $\Xi^+_n$  the cell  of $T$
    on the right of $X_n$. Let  $\lambda^+$ denote the intensity of
    the marked point process $\{(X_n,\Xi_n^+)\}$,
     and $\Xi^{+*}$ its typical mark. Then, by the definition of the
     Palm mark distribution (see e.g.
     Section~\ref{subsec.pal.dis}), we get that
    \begin{eqnarray*}
        \E\, c(X_0^*,X_1^*) &\le& \E\, \nu_1(\partial \Xi^{+*}) \;=
\;        \frac{1}{\lambda^+}\;\E\sum_{X_i\in T^{(1)}\cap [0,1)}\nu_1(\partial \Xi_i^+)\\
        &=&\frac{1}{\lambda^+}\;\E\sum_{\Xi_i\in T}
        \ind_{\{\partial^+ \Xi_i \cap [0,1)\not=\emptyset\}}\nu_1(\partial \Xi_i) \,,
    \end{eqnarray*}
    where $\partial^+\Xi$ denotes that part of the boundary of $\Xi$ with outer unit normal vector in $[\pi/2,3\pi/2)$.
    Thus, applying Campbell's theorem to the latter expression, we have
    \begin{eqnarray*}
        \E c(X_0^*,X_1^*) &\le& \frac{\lambda_T}{\lambda^+}\;\E\nu_1(\partial \Xi^*)
        \int_{\R^2} \ind_{\{\partial^+ \Xi^* + x \cap [0,1)\not=\emptyset\}}\nu_2(dx)\\
        &=&\frac{\lambda_T}{\lambda^+}\;\E\nu_1(\partial \Xi^*)
        \nu_2([0,1)\oplus \partial^+\Xi^*)\,,
    \end{eqnarray*}
    where $\lambda_T=1/\E\,\nu_2(\Xi^*)$. Since $\nu_2([0,1)\oplus
    \partial^+\Xi^*)\le a\, \nu_1(\partial\Xi^*)$ for some
    constant $a<\infty$, this implies that
     $ \E\, c(X_0^*,X_1^*) \le (a\lambda_T/\lambda^+)\,\E\nu_1^2(\partial
     \Xi^*)$.
    Thus, the assertion is shown.

    \vspace{0.5cm} \noindent
 {\bf Ergodicity.} $\;$
    We only prove that $\mathbf{X}$ is ergodic, because the ergodicity of $\mathbf{Y}$ can be shown in the same way.
    Recall that by $\mathcal{I}_\mathcal{S}\subset\mathcal{B}(\mathcal{S})$ we denote the
    $\sigma$-algebra of those subsets of the space $\mathcal{S}$ of double-indexed sequences, which are invariant under the shift
    $\{|X_i^*-X_j^*|\}\longmapsto\{|X_{i+1}^*-X_{j+1}^*|\}$. Furthermore, note  that $\mathbf{X}=h(T_\gamma^{(1)*})$ for some  measurable
    function $h:\mathcal{F}\to \mathcal{S}$, where for any
    tessellation $\tau$ in $\R^2$ and
    $A\in\mathcal{I}_\mathcal{S}$, we have $h(\tau^{(1)})\in A$ if and only if
    $h(\tau^{(1)}-x)\in A$ for all $x\in [0,\infty)$.
    Thus, from the definition of the Palm distribution of the stationary point process
    $\{X_i\}=T^{(1)}\cap \R$ with intensity $2/\pi$, we get for any
    $A\in\mathcal{I}_\mathcal{S}$
    that
    \begin{eqnarray*}
        \mathbb{P}(\mathbf{X}\in A)&=&\mathbb{P}(h(T^{(1)*})\in A)\\
        &=&\frac{\pi}{2}\;\E\sum\limits_{X_i\in T^{(1)}\cap
        B(o,1)\cap L^+} \ind_{A}(h(T^{(1)}-X_i))\\
        &=&\frac{\pi}{2}\;\E\bigl(\ind_{A}(h(T^{(1)}))\;\#\{X_i\in T^{(1)}\cap
        B(o,1)\cap L^+\}\bigr)\\
        &=&\frac{\pi}{2}\;\E\bigl(\ind_{h^{-1}(A)}(T^{(1)})\;\#\{X_i\in T^{(1)}\cap
        B(o,1)\cap L^+\}\bigr)\,.
    \end{eqnarray*}
    On the other hand, since $T_1$ is mixing and $h^{-1}(A)=h^{-1}(A)+x$ for any
    $A\in\mathcal{I}_\mathcal{S}$ and
     $x\in L^+$, we have
    \begin{eqnarray*}
        \mathbb{P}\big(T^{(1)}\in h^{-1}(A)\big)&=&\lim_{|x|\rightarrow \infty,x\in L^+}
        \mathbb{P}\big(T^{(1)}\in h^{-1}(A),T^{(1)}-x\in h^{-1}(A)\big)\\
        &=&\mathbb{P}\big(T^{(1)}\in h^{-1}(A)\big)^2\,,
    \end{eqnarray*}
    which implies that $\mathbb{P}\big(T^{(1)}\!\!\in h^{-1}(A)\big)=0$ or $\mathbb{P}\big(T^{(1)}\!\!\in h^{-1}(A)\big)=1$.
    Thus, altogether, we have
    \begin{eqnarray*}
        \mathbb{P}(\mathbf{X}\in A)
        &=&\mathbb{P}\big(T^{(1)}\in h^{-1}(A)\big)\;\frac{\pi}{2}\;\E\;\#\{X_i\in T^{(1)}\cap
        B(o,1)\cap L^+\} \\
        &=&\mathbb{P}\big(T^{(1)}\in h^{-1}(A)\big)
    \end{eqnarray*}
and, consequently,  $\mathbb{P}(\mathbf{X}\in A)=0$ or
$\mathbb{P}(\mathbf{X}\in A)=1$ for any
    $A\in\mathcal{I}_\mathcal{S}$, which means that  $\mathbf{X}$ is ergodic.

    \vspace{0.5cm} \noindent
 {\bf Uniform integrability.} $\;$
Finally, we show that the family
$\{\widetilde{g}_\gamma,\,\gamma>0\}$  considered in
(\ref{alm.shu.con}) is uniformly $(\nu_1\otimes
\mathbb{P}^*)$-integrable.  From the ergodic theorem for
stationary marked point processes (\cite{Daley0307}, Theorem
12.2.IV), we get that
    \begin{eqnarray*}
        \lim_{\gamma\rightarrow \infty}\frac{1}{\gamma}
        \sum_{\substack{X_i\in T^{(1)}\cap L^+:\\|X_i|\le (r+1)\gamma}}\,\frac{1}{\sin\alpha_i}&=&
        (r+1)\lim_{\gamma\rightarrow \infty}\frac{1}{(r+1)\gamma}
        \sum_{\substack{X_i\in T^{(1)}\cap L^+:\\|X_i|\le (r+1)\gamma}}\,\frac{1}{\sin\alpha_i}\\
        &=&(r+1)\;\E(\sin\alpha^*)^{-1}\,
    \end{eqnarray*}
    almost surely and in $L^1$ since the point process $T^{(1)}\cap \R$ marked with the intersection angles
    is ergodic, which can be shown in the same way as the ergodicity of $\mathbf{X}$. Here $\alpha^*$ denotes the typical
    intersection angle which is distributed according to the density $f_{\alpha^*}(\alpha)=\sin(\alpha)/2$ for $0\le \alpha< \pi$,
    see e.g. \cite{Stoyan95}, p. 288. This yields $\E(\sin\alpha^*)^{-1}=\pi/2<\infty$. Thus
    \begin{eqnarray*}
        0&=&\lim_{\gamma\rightarrow \infty}\E\;
        \bigg|\;\frac{1}{\gamma}\;
        \sum_{\substack{X_i\in T^{(1)}\cap L^+:\\|X_i|\le (r+1)\gamma}}\;\frac{1}{\sin\alpha_i}\;-\;(r+1)\;\E(\sin\alpha^*)^{-1}\bigg|\\
        &=&\lim_{\gamma\rightarrow \infty}\frac{2}{\pi}\;\E \int_\R \ind_{[0,
        X_1^*]}(x)\;\bigg|\;\frac{1}{\gamma}\;
        \sum_{\substack{X_i\in (T^{(1)*}-x)\cap L^+:\\|X_i|\le (r+1)\gamma}}\;\frac{1}{\sin\alpha_i}\;-\;(r+1)\;\E(\sin\alpha^*)^{-1}\bigg|\,dx\,,
    \end{eqnarray*}
   where in the last equality we used the inversion formula for Palm distributions of stationary marked point processes on $\R$;
   see Proposition 11.3 ~(iii) in \cite{Kal02}. In other words, we showed that
    \begin{eqnarray*}
        \ind_{[0, X_1^*]}(x)\;\frac{1}{\gamma}
        \sum_{\substack{X_i\in (T^{(1)*}-x)\cap L^+:\\|X_i|\le
        (r+1)\gamma}}\,\frac{1}{\sin\alpha_i}\;\longrightarrow\;
        (r+1)\ind_{[0, X_1^*]}(x)\;\E(\sin\alpha^*)^{-1}
    \end{eqnarray*}
    in $L^1(\nu_1\otimes \mathbb{P}^*)$ as $\gamma\to\infty$. This
    means in particular that the family $\{h_\gamma,\,\gamma>0\}$
    with
    \[
    h_\gamma(x, T^{(1)*})=   \ind_{[0, X_1^*]}(x)\;\frac{1}{\gamma}
        \sum_{\substack{X_i\in (T^{(1)*}-x)\cap L^+:\\|X_i|\le
        (r+1)\gamma}}\,\frac{1}{\sin\alpha_i}
\]
 is uniformly $(\nu_1\otimes\mathbb{P}^*)$-integrable; see
     Theorem~\ref{lem.uni.int}. Furthermore, we have that
    \begin{eqnarray*}
        \ind_{[0, X_1^*]}(x)g_\gamma(T^{(1)*}-x)\le
        \ind_{[0, X_1^*]}(x)\;\frac{1}{\gamma}\sum_{\substack{X_i\in (T^{(1)*}-x)\cap L^+:\\ |X_i|\le
        (r+1)\gamma}}\,\frac{1}{\sin\alpha_i}\,.
    \end{eqnarray*}
    Thus Lemma~\ref{lem.dom.int} yields that the family
$\{\widetilde{g}_\gamma,\,\gamma>0\}$ considered in
(\ref{alm.shu.con}) is uniformly
 $(\nu_1\otimes\mathbb{P}^*)$-integrable.

\end{appendix}

\acknowledgement This research was supported by Orange Labs through research grant No. $46\,14\,37\,14$.
\endacknowledgement

\bibliographystyle{apt}


\end{document}